\documentclass[12pt,reqno]{amsart}
\usepackage{amsmath,amsfonts,amsthm,amsopn,amssymb}
\usepackage{cite,marginnote}
\usepackage{color,enumitem,graphicx}
\usepackage[colorlinks=true,urlcolor=blue,
citecolor=red,linkcolor=blue,linktocpage,pdfpagelabels,
bookmarksnumbered,bookmarksopen]{hyperref}
\usepackage[english]{babel}

\usepackage[left=2.9cm,right=2.9cm,top=2.8cm,bottom=2.8cm]{geometry}

\numberwithin{equation}{section}

\makeindex

\newtheorem{theorem}{Theorem}[section]

\newtheorem{corollary}{Corollary}[section]
\newtheorem{remark}{Remark}[section]

\newtheorem{lemma}{Lemma}[section]
\newtheorem{iteration lemma}{iteration Lemma}[section]

\newcommand{\s}{\section}

\newcommand{\R}{\mathbb R}

\newcommand{\lab}{\label}
\newcommand{\bt}{\begin{theorem}}
\newcommand{\et}{\end{theorem}}
\newcommand{\bl}{\begin{lemma}}
\newcommand{\el}{\end{lemma}}
\newcommand{\bd}{\begin{definition}}
\newcommand{\ed}{\end{definition}}
\newcommand{\bc}{\begin{corollary}}
\newcommand{\ec}{\end{corollary}}
\newcommand{\bp}{\begin{proof}}
\newcommand{\ep}{\end{proof}}
\newcommand{\bx}{\begin{example}}
\newcommand{\ex}{\end{example}}
\newcommand{\bi}{\begin{exercise}}
\newcommand{\ei}{\end{exercise}}
\newcommand{\bo}{\begin{proposition}}
\newcommand{\eo}{\end{proposition}}
\newcommand{\br}{\begin{remark}}
\newcommand{\er}{\end{remark}}
\newcommand{\beq}{\begin{equation}}
\newcommand{\eeq}{\end{equation}}
\newcommand{\ba}{\begin{align}}
\newcommand{\ea}{\end{align}}
\newcommand{\bn}{\begin{enumerate}}
\newcommand{\en}{\end{enumerate}}
\newcommand{\bg}{\begin{align*}}
\newcommand{\bcs}{\begin{cases}}
\newcommand{\ecs}{\end{cases}}

\newcommand{\bean}{\begin{eqnarray*}}
\newcommand{\eean}{\end{eqnarray*}}

\def\N{\mathbb{N}}

\def\R{\mathbb{R}}

\def\bd{\mathrm{bd}\,}

\title[Normalized solutions to Kirchhoff type equations]{Normalized solutions to Kirchhoff type equations with a critical growth nonlinearity}

\author[J. Zhang]{Jian Zhang}
\author[J. J. Zhang]{Jianjun Zhang}
\author[X.~X.~Zhong]{Xuexiu Zhong}

\address[J. \ Zhang]{\newline\indent College of Science, China University of Petroleum
\newline\indent
Qingdao 266580, Shandong, PR China}
\email{\href{mailto:zjianmath@163.com}{zjianmath@163.com}}

\address[J. J. \ Zhang]{\newline\indent College of Mathematics and Statistics, Chongqing Jiaotong University
\newline\indent
Chongqing 400074, PR China}
\email{\href{mailto:zhangjianjun09@tsinghua.org.cn}{zhangjianjun09@tsinghua.org.cn}}

\address[X.~X.~Zhong]{\newline\indent South China Research Center for Applied Mathematics and Interdisciplinary Studies
\newline\indent
South China Normal University
\newline\indent
Guangzhou 510631, PR China}
\email{\href{mailto:zhongxuexiu1989@163.com}{zhongxuexiu1989@163.com}}

\thanks{Xuexiu Zhong was supported by the NSFC (No.12271184), Guangdong Basic and Applied Basic Research Foundation (2021A1515010034),Guangzhou Basic and Applied Basic Research Foundation(No.202102020225). Jianjun Zhang was supported by NSFC (No.11871123).}

\date{}

\begin{document}

\begin{abstract}
 In this paper, we are concerned with normalized solutions of the Kirchhoff type equation
 \begin{equation*}
-M\left(\int_{\R^N}|\nabla u|^2\mathrm{d} x\right)\Delta u  = \lambda u +f(u) \ \ \mathrm{in} \ \ \mathbb{R}^N
\end{equation*}
with $u \in S_c:=\left\{u \in H^1(\R^N): \int_{\R^N}u^2 \mathrm{d}x=c^2\right\}$. When $N=2$ and $f$ has exponential critical growth at infinity, normalized mountain pass type solutions are obtained via the variational methods. When $N \ge 4$, $M(t)=a+bt$ with $a$, $b>0$ and $f$ has Sobolev critical growth at infinity, we investigate the existence of normalized ground state solutions, mountain pass type solutions and local constraint minimizer with positive energy. Moreover, the non-existence of normalized solutions is also considered.
\end{abstract}

\maketitle

\s{Introduction and results}
\renewcommand{\theequation}{1.\arabic{equation}}

\subsection{Background} The Kirchhoff type problem appears as models of several physical phenomena. For example, it
is related to the stationary analogue of the equation:
\begin{equation}\label{1-1}
\rho \frac{\partial^2 u}{\partial t^2}-\left(\frac{P_0}{h}+\frac{E}{2L}\int_0^L\left|\frac{\partial u}{\partial x}\right|^2\mathrm{d}x\right)\frac{\partial^2 u}{\partial x^2}=0,
\end{equation}
where $u$ is the lateral displacement at $x$ and $t$, $E$ is the Young modulus, $\rho$ is the mass density, $h$ is the cross-section area, $L$ is the length, $P_0$ is the initial axial tension.
For more background, see \cite{ABG,K} and the references
therein. Because of the presence of the nonlocal term, the Kirchhoff type equation is no longer a pointwise identity, which causes additional mathematical difficulties.

In this paper, we study solutions to the following Kirchhoff type equation with critical growth nonlinearities:
\begin{equation}\label{1-2}
-M\left(\int_{\R^N}|\nabla u|^2\mathrm{d} x\right)\Delta u =\lambda u + f(u) \ \ \mathrm{in} \ \ \mathbb{R}^N
\end{equation}
having prescribed mass
\begin{align}\lab{eq:constraint}
\int_{\R^N}u^2 \mathrm{d}x=c^2,
\end{align}
where $N=2$ or $N \ge 4$, $\lambda \in \R$ is an unknown parameter appearing as a Lagrange multiplier and $f$ is of critical growth.
For the case of $N=2$, we consider the general Kirchhoff type case. While for $N\geq 4$, we consider the classical Kirchhoff case that $M(t)=a+bt$ with $a,b>0$.

A solution $u$ to Eq.\eqref{1-2} satisfying \eqref{eq:constraint} is called a normalized solution. One can search for normalized solutions of \eqref{1-2}  by
studying critical points of the functional
\beq\lab{eq:def-energy-functional}
I(u)=\frac{1}{2}\widehat{M}\left(\int_{\R^N}|\nabla u|^2\mathrm{d} x\right)-\int_{\R^N}F(u)\mathrm{d} x,
\eeq
constrained on the $S_c$,
where $\widehat{M}(t):=\int_0^t M(s)\mathrm{d}s$, $F(u):=\int_0^u f(s)\mathrm{d}s$ and
\beq\lab{eq:def-Sc}
S_c:=\left\{u \in H^1(\R^N):\int_{\R^N}u^2 \mathrm{d}x=c^2\right\}.
\eeq

For $c>0$, if $I|_{S_c}'(u_c)=0$, then there exists some $\lambda_c \in \R$ such that $(u_c,\lambda_c)$ is a couple of solution to (\ref{1-2}), where $\lambda_c$ is the associated Lagrange multiplier. For any $c>0$, define
\beq\lab{def:Ic}
I_c:=\inf_{u \in S_c}I(u).
\eeq
We say that $u_0$ is a normalized ground state solution of \eqref{1-2} if
\begin{align*}
I(u_0)=\inf\left\{I(u): I|_{S_c}'(u)=0, u \in S_c\right\}.
\end{align*}
In particular, if $u_0$ attains $I_c$, then $u_0$ is a normalized ground state solution of \eqref{1-2}.

The normalized solution has important physical background and
it has attracted attentions of many researchers in recent years. For the convenience, we refer the readers to the following papers, and the
references therein, \cite{AW,B,BJ,BJS,BS,BZZ,BJ,BJL,BCGJ,BES,GJ,GZ,GLW,GS,GWZZ,HT,J,JL,JeL,JLW,JZZ,LZ,NTV,PV,S1,S2,WW}. In particular, L. Jeanjean \cite{J} exploited the mountain pass lemma to study normalized solutions of Schr\"{o}dinger equations for the $L^2$-supercritical case. N. Soave \cite{S2} considered normalized ground state solutions of Schr\"{o}dinger equations with combined nonlinearities for the Sobolev critical case. Jeanjean and Lu \cite{JeL}
investigated normalized local minimizers and mountain pass type solutions for a coercive problem. Jeanjean, Zhang and Zhong \cite{JZZ} introduced a global branch approach to study normalized solutions of the Schr\"{o}dinger equation.
There are relatively few results about normalized solutions of Kirchhoff type equations.
In \cite{Ye,Y}, the author considered normalized solutions of
\eqref{1-2} for the case of $N \leq 3$, $M(t)=a+bt$ with $a$, $b>0$ and $f(u)=|u|^{p-2}u$ with $p \in (2,2^*)$.
Precisely, for $p \in (2,\frac{2N+8}{N})$, the authors studied the existence and non-existence of global constraint minimizers; for $p \in (\frac{2N+4}{N},\frac{2N+8}{N})$, the authors studied local constraint minimizers;
for $p \in [\frac{2N+8}{N},2^*)$, the authors studied mountain pass type solutions.
In \cite{ZZ}, the authors improved the results in \cite{Ye}. In \cite{QZ}, the authors investigated
the exact number and expressions of positive normalized solutions for Kirchhoff equations.
In \cite{L,H}, the authors studied normalized solutions of Kirchhoff equations with critical growth nonlinearity in dimension three.

To the best of our knowledge, there are no results on normalized solutions of Kirchhoff type equations with exponential critical growth in dimension two.
Here we recall that the nonlinearity $f$ has exponential subcritical growth if for any $\alpha>0$,
\begin{align*}
\lim_{u \rightarrow +\infty}\frac{f(u)}{e^{\alpha u^2}}=0
\end{align*}
and the nonlinearity $f$ has exponential critical growth if there exists $\alpha_0>0$ such that
\begin{align*}
\lim_{u \rightarrow +\infty}\frac{f(u)}{e^{\alpha u^2}}=
\left\{
\begin{aligned}
 &0,       &\forall\alpha>\alpha_0,\\
 &+\infty, &\forall\alpha<\alpha_0.\\
\end{aligned}
\right.
\end{align*}
We notice that in \cite{AJM}, the authors used the idea due to L. Jeanjean \cite{J} to study normalized solutions of Schr\"{o}dinger equations (i.e., $M(t)\equiv 1$ in \eqref{1-2}) for $N=2$, where $f$ is of exponential critical growth satisfying the following technical condition.
\begin{itemize}
\item [$(f)$] There exists $p>4$ and $\mu>0$ such that $\mbox{sgn}(t) f(t) \ge \mu |t|^{p-1}$ for $t \ne 0$, where $\mbox{sgn}(t)=1$ if $t>0$ and $\mbox{sgn}(t)=-1$ if $t<0$.
\end{itemize}
We remark that in \cite{AJM}, the condition $(f)$ with $\mu>0$ being large is essential, which enables the  upper bound on the energy being controlled. It is natural to ask if the restriction can be removed. This is one motivation of the present paper.

We recall that in \cite{dMR},
the authors studied positive solutions of the Dirichlet problem:
\begin{equation*}
-\Delta u = f(x,u) \ \ \mathrm{in} \ \
\Omega, \ \ u=0 \ \mathrm{on} \ \partial \Omega.
\end{equation*}
In this classical paper, to estimate the upper bound on the energy, the author introduced the following more natural condition:
\begin{itemize}
\item [$(f')$] There exists $\beta> \frac{4}{3 \alpha_0 d^2}$ such that
\begin{align*}
\lim_{u \rightarrow +\infty}\frac{f(x,u)u}{e^{\alpha_0 u^2}} \ge \beta,
\end{align*}
where $d$ is the radius of the largest open ball in $\Omega$.
\end{itemize}
This argument is based on the Moser sequence of functions and the proof by the contradiction.
Similar results can be found in \cite{FS,RS}. However, the arguments
cannot be applied directly to study normalized solutions. This is one of the difficulties. We give a direct argument, which enables us to
get the desired upper bound estimate on the energy using the condition $(f_3)$.

Another difficulty is the presence of the nonlocal term, which brings an additional obstacle in proving the $PSP$-condition in Lemma \ref{lemma:20220928-hbl4}.
It is well known that the Trudinger-Moser inequality is essential to deal with the loss of the compactness caused by the exponential critical nonlinearity. When using this inequality, a key step is to give a suitable uniformly upper bound on the $H^1$-norm of the sequence.
However, the nonlocal term prevents us from using the upper bound on the energy and the Ambrosetti-Rabinowitz type condition to deduce the suitable $H^1$-norm estimate. We firstly give a compactness result in Lemma \ref{lemma:20221217}, based on which we obtain the desired $H^1$-norm estimate.
Then we use the Trudinger-Moser inequality to establish a Br\'{e}zis-Lieb type result in Lemma \ref{cro:20221018-c1}
and solve the compactness of the Kirchhoff type problem through
a brief argument. We believe that the Br\'{e}zis-Lieb type result can also be used to study related compact problems in a non-radial setting.

\subsection{Main results}

Now we state our result in dimension two.
We consider the so-called general Kirchhoff type equation (\ref{1-2})-(\ref{eq:constraint}).
For the nonlocal term $M(t)$, we assume the following conditions:
\begin{itemize}
\item [$(M_1)$] $M \in C([0,+\infty))\cap C^1([R,+\infty))$ for some $R>0$, $M_0:=M(0)>0$ and $M(t)$  increases respect to $t \in \R^+$.
\item [$(M_2)$] There exists some $\theta>0$ such that $\displaystyle\widehat{M}(t)-\frac{1}{\theta+1}M(t)t \ge 0$ for $t \in \R^+$ and $\displaystyle\lim_{t \rightarrow +\infty}\left[\widehat{M}(t)-\frac{1}{\theta+1}M(t)t\right]=+\infty$, where $\widehat{M}(t)=\int_0^t M(s) \mathrm{d}s$.
\end{itemize}
For the nonlinear term $f(u)$, we suppose the following:
\begin{itemize}
\item [$(f_1)$] $f \in C(\R^+, \R^+)$ and $\displaystyle\lim_{u \rightarrow 0+}\frac{f(u)}{u^3}=0$.
\item [$(f_2)$] There exists $\alpha_0>0$ such that
\begin{align*}
\lim_{u \rightarrow +\infty}\frac{f(u)}{e^{\alpha u^2}}=
\left\{
\begin{aligned}
 &0,       &\forall\alpha>\alpha_0,\\
 &+\infty, &\forall\alpha<\alpha_0.\\
\end{aligned}
\right.
\end{align*}
\item [$(f_3)$] There exists $\beta>0$ such that
$\displaystyle\liminf_{u \rightarrow +\infty}\frac{f(u)u}{e^{\alpha_0 u^2}} \geq \beta$.
\item [$(f_4)$] $\frac{1}{2 \theta +4}f(u)u-F(u) \ge 0$ for all $u \in \R^+$, where $F(u)=\int_0^u f(s) \mathrm{d}s$ and $\theta>0$ is given by $(M_2)$.
 \item [$(f_5)$] $\limsup_{u \rightarrow +\infty}\frac{F(u)}{f(u)u}=0$.
\end{itemize}

\bt\lab{Theorem 1.1} Assume that $(M_1)$-$(M_2)$ and $(f_1)$-$(f_5)$ hold.
Then problem (\ref{1-2})-(\ref{eq:constraint}) has a couple of solution $(u_c,\lambda_c) \in H^1(\R^2) \times \R$
such that $\lambda_c<0$ and $u_c$ is a mountain pass type critical point of $I|_{S_c}$.
\et

\br We recall that in \cite{dMR}, the compactness result in Lemma $2.1$ plays an important role in proving the $PS$-condition of the Dirichlet problem with exponential critical growth in $\R^2$. To get the compactness result, the authors introduced the following classical condition:
\begin{itemize}
\item [$(f_5')$] There exist $u_0$, $L_0>0$ such that $F(u) \le L_0 f(u)$ for $u \ge u_0$.
\end{itemize}
To the best of our knowledge, the condition $(f_5')$ was later used to study all critical problems in $\R^2$. A natural question is if $(f_5')$ can be weaken or removed? In this paper, we introduce a new condition $(f_5)$, which can replace $(f_5')$ to study the compactness of the problem with exponential critical growth in $\R^2$.
\er

For the case of $N \ge 4$, there are few results about normalized solutions to Kirchhoff type equations. We notice that in \cite{Z},
the authors studied normalized solutions of
(\ref{1-2})-(\ref{eq:constraint}),
where $M(t)=a+bt$ with $a$, $b>0$ and $f$ satisfies the following conditions:
\begin{itemize}
\item [$(f'_1)$] $f \in C^1(\R^+, \R^+)$ and $f(u)>0$ for $u>0$.
\item [$(f'_2)$] There exist $2<\alpha$, $\beta<2^*:=\frac{2N}{N-2}$ such that $\displaystyle\lim_{u \rightarrow 0+}\frac{f'(u)}{u^{\alpha-2}}=\mu_1(\alpha-1)>0$ and $\displaystyle\lim_{u \rightarrow +\infty}\frac{f'(u)}{u^{\beta-2}}=\mu_2(\beta-1)>0$.
\item [$(f'_3)$] $-\Delta u=f(u)$ has no positive radial decreasing classical solution in $\R^N$.
\end{itemize}
Let $U$ be the unique positive solutions of the following equation:
\begin{align*}
-\Delta U + U=\mu_1 U^{\alpha-1} \ \ \mathrm{in} \ \R^N, \ \ \lim_{|x| \rightarrow \infty}U(x)=0.
\end{align*}
By using the global branch approach developed by Jeanjean et al. in \cite{JZZ} and the  Azzollini's correspondence in \cite{A}, which provided a homeomorphism between ground state
solutions of autonomous Kirchhoff type equations and related local semilinear
elliptic equations, the authors obtained the following results in \cite{Z}.

\bt\lab{Theorem 1} (The case $N=4$) Let $M(t)=a+bt$ with $a$, $b>0$ and $(f'_1)$-$(f'_3)$ hold.
\begin{itemize}
\item[(1)] If $\alpha \in (3,4)$, then there exists $c^*>0$ such that for any $c>c^*$, problem \eqref{1-2}-\eqref{eq:constraint} has at least two distinct positive normalized solutions $(\lambda_i,u_{\lambda_i}) \in (-\infty,0)\times H_r^1(\R^4)$, $i=1$, $2$.
\item[(2)]If $\alpha =3$, then for any $c > a \|U\|_2$, problem \eqref{1-2}-\eqref{eq:constraint} has at least one positive normalized solution $(\lambda,u_\lambda) \in (-\infty,0)\times H_r^1(\R^4)$.
\item[(3)] If $\alpha \in (2,3)$, then for any $c>0$, problem \eqref{1-2}-\eqref{eq:constraint} has at least one positive normalized solution $(\lambda,u_\lambda) \in (-\infty,0)\times H_r^1(\R^4)$.
\end{itemize}
 \et

\bt\lab{Theorem 2} (The case $N \ge 5$) Let $M(t)=a+bt$ with $a$, $b>0$ and $(f'_1)$-$(f'_3)$ hold.
\begin{itemize}
\item[(1)] If $\alpha \in (2+\frac{4}{N},2^*)$, then there exists $c^*>0$ such that for any $c>c^*$, problem \eqref{1-2}-\eqref{eq:constraint} has at least two distinct positive normalized solutions $(\lambda_i,u_{\lambda_i}) \in (-\infty,0)\times H_r^1(\R^N)$, $i=1$, $2$.
\item[(2)] If $\alpha =2+\frac{4}{N}$, then for any $c > a^{\frac{N}{4}} \|U\|_2$, problem \eqref{1-2}-\eqref{eq:constraint} has at least one positive normalized solution $(\lambda,u_\lambda) \in (-\infty,0)\times H_r^1(\R^N)$.
\item[(3)]  If $\alpha \in (2,2+\frac{4}{N})$, then for any $c>0$, problem \eqref{1-2}-\eqref{eq:constraint} has at least one positive normalized solution $(\lambda,u_\lambda) \in (-\infty,0)\times H_r^1(\R^N)$.
\end{itemize}
 \et

\bt\lab{Theorem 3}  Let $M(t)=a+bt$ with $a$, $b>0$, $N \ge 4$ and $(f'_1)$-$(f'_3)$ hold. If $\alpha \in [2+\frac{4}{N},2^*)$, then there exists $c_0>0$ such that problem \eqref{1-2}-\eqref{eq:constraint} has no positive solution for $c<c_0$.
 \et

For the critical case, we cannot use the Azzollini's correspondence in \cite{A} to study the normalized solution problem, since the mass is not preserved via the transformation.
And we also can not use the global branch approach in \cite{JZZ} since the presence of the Sobolev critical term. Motivated by the above facts,
we consider the following Kirchhoff type equation with combined nonlinearities
 \begin{equation}\label{1.5}
-\left(a+b\int_{\R^N}|\nabla u|^2\mathrm{d} x\right)\Delta u =\lambda u + |u|^{p-2}u+ |u|^{2^*-2}u \ \ \mathrm{in} \ \ \mathbb{R}^N,
\end{equation}
for some $\lambda \in \R$ having prescribed  mass
\begin{align}\label{1.6}
\int_{\R^2}u^2 \mathrm{d}x=c^2,
\end{align}
in the Sobolev critical case, where $N \ge 4$, $a$, $b>0$ and $2<p<2^*:=\frac{2N}{N-2}$.

The Sobolev critical term complicates the problem. As we will see, the interplay between Sobolev subcritical and Sobolev critical nonlinearities strongly affects the structure of the functional, the compactness of the problem and the existence of normalized solutions. For the case of $N \ge 4$,
since $2^*:=\frac{2N}{N-2} \le 4$, the critical term will be affected by the nonlocal term.
We give a technical Lemma \ref{lemma:20220923-hhl1},
which helps to deal with the structure of the functional and solve the loss of the compactness caused by the critical term.

The $L^2$-critical exponent $2+\frac{4}{N}$ plays an important role in studying normalized solutions.
For the critical case, the nonlinearity involves a critical term $|u|^{2^*-2}u$ growing faster than $|u|^{\frac{4}{N}}u$ at infinity. Furthermore, for $N \ge 4$, we have that $2^* \le 2+\frac{8}{N}$, where $2+\frac{8}{N}$ is the so called $L^2$-critical exponent for Kirchhoff equations, which is quite different from the case $N=3$. In this paper, we study the problem for the case of $p>2+\frac{4}{N}$, $p=2+\frac{4}{N}$ and $p<2+\frac{4}{N}$ respectively. Firstly, we study the non-existence of normalized solutions.
Secondly, we study the existence of normalized ground state solutions. Thirdly, we study normalized mountain pass type solutions and normalized local constraint minimizer with positive energy.

\vskip0.13in

We introduce the following Gagliardo-Nirenberg inequality.

\bl\lab{lemma:20220923-l1} (\cite{W}) Let $p \in [2,\frac{2N}{N-2})$ if $N \ge 3$ and $p \ge 2$ if $N=1$, $2$. Then
\begin{align}\label{1.7}
\|u\|_p \le \left(\frac{p}{2 \|Q\|_2^{p-2}}\right)^{\frac{1}{p}}\|\nabla u\|_2^{\frac{N(p-2)}{2p}}\|u\|_2^{1-\frac{N(p-2)}{2p}},
\end{align}
with equality only for $u=Q$, where up to a translation, $Q$ is the unique ground state solution of
\begin{align}\label{1.8}
-\frac{N(p-2)}{4}\Delta Q + \left(1+\frac{(p-2)(2-N)}{4}\right)Q=|Q|^{p-2}Q, \ \ x \in \R^N.
\end{align}
\el

Here comes our main results for $N\geq 5$.

\bt\lab{th:main-zzzt1} Let $N \ge 5$, $\left(\frac{2a}{4-2^*}\right)^{\frac{4-2^*}{2}} \left(\frac{2b}{2^*-2}\right)^{\frac{2^*-2}{2}}>\frac{1}{S^{\frac{2^*}{2}}}$ and $p \in (2,2^*)$.
\begin{itemize}
\item[(i)] If $p \in \left(2+\frac{4}{N},2^*\right)$, then
\begin{itemize}
 \item [(i-1)] there exists $c_1 \in (0,+\infty)$ such that $I_c=0$ is not attained for $c \in (0,c_1)$, $I_c=0$ is attained for $c =c_1$ and $I_c<0$ is attained for $c>c_1$. In particular, for $c \ge c_1$, problem (\ref{1.5})-\eqref{1.6} has a normalized ground state solution.
 \item [(i-2)] for $c \ge c_1$, problem (\ref{1.5})-\eqref{1.6} has a couple of solution $(u_c,\lambda_c) \in H^1(\R^N) \times \R$
such that  $\lambda_c<0$ and $u_c$ is a mountain pass type critical point of $I|_{S_c}$.
\item [(i-3)] there exists $\bar{\eta}>0$ such that for $c \in [c_1-\bar{\eta},c_1)$, problem (\ref{1.5})-\eqref{1.6} has two couples of solutions $\left(u_c^{(i)},\lambda_c^{(i)}\right) \in H^1(\R^N) \times \R$, $i=1$, $2$ such that $\lambda_c^{(1)}<0$ and $u_c^{(1)}$ is a mountain pass type critical point of $I|_{S_c}$;  $\lambda_c^{(2)}<0$ and $u_c^{(2)}$ is a local constraint minimizer of $I|_{S_c}$ with positive energy.
\item [(i-4)] there exists $c_0 \in (0,c_1-\bar{\eta}]$ such that for $c \in (0,c_0)$, problem (\ref{1.5})-\eqref{1.6} has no normalized solution.
In particular, $I|_{S_c}$ has no constraint critical points.
\end{itemize}
\item[(ii)] If $p = 2+\frac{4}{N}$, then
\begin{itemize}
 \item [(ii-1)] there exists $c_1>0$ such that $I_c=0$ is not attained for $c \in (0,c_1)$, $I_c=0$ for $c =c_1$ and $I_c<0$ is attained for $c>c_1$. In particular, for $c > c_1$, problem (\ref{1.5})-\eqref{1.6} has a normalized ground state solution.
\item [(ii-2)] there exists $c_0\geq \|Q\|_2 \left[a -\frac{4-2^*}{2S^{\frac{N}{N-4}}}\left(\frac{2^*-2}{2b}\right)^{\frac{2}{N-4}}\right]^{\frac{N}{4}}$ such that problem (\ref{1.5})-\eqref{1.6} has no normalized solution provided $c\in (0,c_0)$. Consequently, $I|_{S_c}$ has no constraint critical points.
\end{itemize}
\item[(iii)] If $p \in \left(2,2+\frac{4}{N}\right)$, then $I_c<0$ is attained for all $c>0$. In particular, problem (\ref{1.5})-\eqref{1.6} has a normalized ground state solution.
\end{itemize}
 \et

For the case of $N=4$, we have the following results.

\bt\lab{th:main-zzzt2} Let $N=4$, $b>\frac{1}{S^2}$ and $p \in (2,4)$.\\
\begin{itemize}
\item[(i)]If $p \in (3,4)$, then
\begin{itemize}
\item [(i-1)] there exists $c_1 \in (0,+\infty)$ such that $I_c=0$ is not attained for $c \in (0,c_1)$, $I_c=0$ is attained for $c =c_1$ and $I_c<0$ is attained for $c>c_1$. In particular, for $c \ge c_1$, problem (\ref{1.5})-\eqref{1.6} has a normalized ground state solution.
\item [(i-2)] for $c \ge c_1$, problem (\ref{1.5})-\eqref{1.6} has a couple of solution $(u_c,\lambda_c) \in H^1(\R^4) \times \R$
such that $\lambda_c<0$ and $u_c$ is a mountain pass type critical point of $I|_{S_c}$.
\item [(i-3)] there exists $\bar{\eta}>0$ such that for $c \in [c_1-\bar{\eta},c_1)$, problem (\ref{1.5})-\eqref{1.6} has two couples of solutions $\left(u_c^{(i)},\lambda_c^{(i)}\right) \in H^1(\R^4) \times \R$, $i=1$, $2$ such that $\lambda_c^{(1)}<0$ and $u_c^{(1)}$ is a mountain pass type critical point of $I|_{S_c}$;  $\lambda_c^{(2)}<0$ and $u_c^{(2)}$ is a local constraint minimizer of $I|_{S_c}$ with positive energy.
\item [(i-4)] there exists some $c_0\geq \frac{a \|Q\|_2^{\frac{p-2}{4-p}}}{(4-p)(p-2)^{\frac{1}{4-p}}}\left[\frac{1}{p-3}\left(b-\frac{1}{S^2}\right)\right]^{\frac{p-3}{4-p}}$ such that  problem (\ref{1.5})-\eqref{1.6} has no normalized solution provided $c \in \left(0,c_0\right)$.
In particular, $I|_{S_c}$ has no constraint critical points.
\end{itemize}
\item[(ii)] If $p =3$, then
\begin{itemize}
\item [(ii-1)] $I_c=0$ is not attained for $c \in (0,a\|Q\|_2]$ and $I_c<0$ is attained for $c>a\|Q\|_2$. In particular, for $c>a\|Q\|_2$, problem (\ref{1.5})-\eqref{1.6} has a normalized ground state solution.
\item [(ii-2)] for $c \in (0,a\|Q\|_2]$, problem (\ref{1.5})-\eqref{1.6} has no normalized solution. In particular, $I|_{S_c}$ has no constraint critical points.
\end{itemize}
\item[(iii)] If $p \in (2,3)$, then $I_c<0$ is attained for all $c>0$. In particular, problem (\ref{1.5})-\eqref{1.6} has a normalized ground state solution.
\end{itemize}
 \et

\br Compared with \cite{Z}, we establish a more specific relationship between the range of $c$ and the existence, non-existence, multiplicity of normalized solutions. Furthermore, we get types of normalized solutions.
\er

The outline of this paper is as follows: in Section \ref{sec:caseN=2}, we
study the case $N=2$; in Section \ref{sec:caseNgeq5}, we study the case $N \ge 5$;
in Section \ref{sec:caseN=4}, we study the case $N =4$.
\vskip0.13in

\noindent{\bf Notations:}
\begin{itemize}
\item [$\bullet$] Denote $\|u\|_s:=\big(\int_{\mathbb{R}^N}|u|^s\mathrm{d}x\big)^{\frac{1}{s}}$, where $1 \le
s < \infty$.
\item [$\bullet$] Denote $H^1(\R^N)$ the Hilbert space with the norm
\begin{align*}
\|u\|_{H^1}:=\left(\|\nabla u\|_2^2+\|u\|_2^2\right)^{\frac{1}{2}}.
\end{align*}
\item [$\bullet$] Denote $D^{1,2}(\mathbb{R}^N)$ the Sobolev space with the norm
\begin{align*}
\|u\|_{D^{1,2}}:=\|\nabla u\|_2.
\end{align*}
\item [$\bullet$] Denote $C$ a universal positive constant (possibly different).
\end{itemize}

\s{The case $N=2$}\lab{sec:caseN=2}

\renewcommand{\theequation}{2.\arabic{equation}}

Let
\begin{align*}
H_r^1(\R^N):=\{u \in H^1(\R^N): u(x)=u(|x|)\}.
\end{align*}
It is well known that the embedding $H_r^1(\mathbb{R}^N) \hookrightarrow L^t(\R^N)$ is compact for all $t \in (2,2^*)$, where $2^*=\infty$ if $N=2$ and $2^*=\frac{2N}{N-2}$ if $N \ge 3$. To deal with the compactness of the problem, instead of $H^1(\R^N)$, we work in the subspace $H_r^1(\R^N)$.

Without loss of generality, we assume that $f(u)=0$ for $u \le 0$. Define the functional $I: H_r^1(\R^N) \rightarrow \R$ as follows:
\begin{align*}
I(u)=\frac{1}{2}\widehat{M}(\|\nabla u\|_2^2)-\int_{\R^N}F(u) \mathrm{d}x.
\end{align*}
Let $S_{r,c}:=S_c \cap H_r^1(\R^N)$. Define
\begin{align*}
P_{r,c}:=\left\{u \in S_{r,c}:G(u)=0\right\},
\end{align*}
where
\begin{align}\lab{eq:20220923-ze1}
G(u)=M(\|\nabla u\|_2^2)\|\nabla u\|_2^2+N\int_{\R^N}F(u) \mathrm{d}x-\frac{N}{2}\int_{\R^N}f(u)u \mathrm{d}x.
\end{align}
Define $H:=H^1(\R^N) \times \R$ with the scalar product
\begin{align*}
\langle\cdot,\cdot\rangle_H=\langle\cdot,\cdot\rangle_{H^1(\R^N)}+\langle\cdot,\cdot\rangle_{\R}.
\end{align*}
Then the norm is
\begin{align*}
\|\cdot\|=\left(\|\cdot\|_{H^1}^2+|\cdot|_{\R}^2\right)^{\frac{1}{2}}.
\end{align*}
Define $T: H \rightarrow H^1(\R^N)$ as follows:
\begin{align*}
T(u,s)=e^{\frac{Ns}{2}}u(e^s x).
\end{align*}
Then $\|T(u,s)\|_2=\|u\|_2$ for all $u \in H^1(\R^N)$ and $s \in \R$. \\

We introduce the following Trudinger-Moser inequality.

\bl\lab{lemma:20220928-hbl1} (\cite{d,R,RS}) If $u \in H^1(\R^2)$ and $\alpha>0$, then
\begin{align}\label{2-1}
\int_{\R^2}\left(e^{\alpha u^2}-1\right)\mathrm{d}x < \infty.
\end{align}
Moreover, for any fixed $\tau>0$, there exists a constant $C>0$ such that
\begin{align}\label{2-2}
\sup_{u \in H^1(\R^2): \|\nabla u\|_2^2 + \tau \|u\|_2^2 \le 1}\int_{\R^2}\left(e^{4 \pi u^2}-1\right)\mathrm{d}x \le C.
\end{align}
\el

\bl\lab{lemma:20220929} For any $\varepsilon>0$, there exist $C_{1,\varepsilon}$, $C_{2,\varepsilon}>0$ such that
\begin{align*}
M(t) \le C_{1,\varepsilon}+C_{2,\varepsilon} t^{\theta+\varepsilon}, \ \ \forall\  t \in \R^+.
\end{align*}
\el

\bp By $(M_1)$-$(M_2)$, we get
\begin{align*}
\frac{1}{\theta+1} \le \limsup_{t \rightarrow +\infty}\frac{\widehat{M}(t)}{M(t)t}=\limsup_{t \rightarrow +\infty}\frac{M(t)}{M(t)+M'(t)t}.
\end{align*}
Then for any $\varepsilon>0$, there exists $T_\varepsilon>0$ such that for $t \ge T_\varepsilon$,
\begin{align*}
\frac{1}{\varepsilon+\theta+1} \le \frac{M(t)}{M(t)+M'(t)t},
\end{align*}
from which we derive that
\begin{align*}
\frac{\theta+\varepsilon}{t} \ge \frac{M'(t)}{M(t)}, \ \ \forall \ t \ge T_\varepsilon.
\end{align*}
Thus,
\begin{align*}
\int_{T_\varepsilon}^t \frac{\theta+\varepsilon}{s} \mathrm{d}s \ge  \int_{T_\varepsilon}^t \frac{M'(s)}{M(s)}\mathrm{d}s.
\end{align*}
Let $C_\varepsilon=\frac{M(T_\varepsilon)}{T_\varepsilon^{\theta+\varepsilon}}$, we have that
\begin{align*}
M(t) \le C_\varepsilon t^{\theta+\varepsilon}, \ \ \forall t \ge T_\varepsilon.
\end{align*}
Together with $(M_1)$, we get the result.
\ep

\subsection{Mountain pass structure}\lab{subsection:structure}

By the structure of $I$, we get the following results.
\bl\lab{lemma:20220928-hbl5} Let $u \in S_{r,c}$. Then
\begin{itemize}
\item [$(1)$] $\|\nabla T(u,s)\|_2 \rightarrow 0$ and $I(T(u,s)) \rightarrow 0$ as $s \rightarrow -\infty$;
\item [$(2)$] $\|\nabla T(u,s)\|_2 \rightarrow +\infty$ and $I(T(u,s)) \rightarrow -\infty$ as $s \rightarrow +\infty$.
\end{itemize}
\el

\bp Obviously, $\lim_{s \rightarrow -\infty}\|\nabla T(u,s)\|_2 = 0$
and $\lim_{s \rightarrow -\infty}\|T(u,s)\|_p=0$ for any $p>2$. Let $\alpha>\alpha_0$ and $q>4$. By $(f_1)$-$(f_2)$, for any $\varepsilon>0$, there exists $C_\varepsilon>0$ such that
\begin{align}\label{2-5}
|f(u)|\leq\varepsilon |u|^3+C_\varepsilon|u|^{q-1}(e^{\alpha u^2}-1), \ \ \forall \ u \in \R.
\end{align}
Then
\begin{align}\label{2-6}
|F(u)| \leq\frac{\varepsilon}{4}|u|^4+\frac{C_\varepsilon}{q}|u|^{q}(e^{\alpha u^2}-1), \ \ \forall \ u \in \R.
\end{align}
By (\ref{2-6}) and Lemma \ref{lemma:20220928-hbl1}, we obtain that
$\lim_{s \rightarrow -\infty}I(T(u,s)) = 0$.
Let $p>2 \theta+3$. By $(f_1)$-$(f_2)$, there exist $c_1$, $c_2>0$ such that
\begin{align}\label{2-9}
f(u) \ge c_1 u^p -c_2 u^3, \ \ \forall \ u \in \R^+.
\end{align}
We note that $\|\nabla T(u,s)\|_2^2 =e^{2s}\|\nabla u\|_2^2$
and $\|T(u,s)\|_{p+1}^{p+1}=e^{s(p-1)}\|u\|_{p+1}^{p+1}$.
Then by (\ref{2-9}) and Lemma \ref{lemma:20220929}, we get $\lim_{s \rightarrow +\infty}\|\nabla T(u,s)\|_2 = +\infty$ and $\lim_{s \rightarrow +\infty}I(T(u,s))= -\infty$.
\ep

\bl\lab{lemma:20220928-hbl6} There exists $K_c>0$ small such that
\begin{align*}
0< \sup_{u \in A_c}I(u) < \inf_{u \in B_c}I(u),
\end{align*}
where
\begin{align*}
A_c:=\{u \in S_{r,c}: \|\nabla u\|_2^2 \le K_c\}, \ \
B_c:=\{u \in S_{r,c}: \|\nabla u\|_2^2 = 2 K_c\}.
\end{align*}
Moreover, $I(u)>0$ for all $u \in A_c$.
\el

\bp  We note that
\begin{align}\label{2.49}
I(v)-I(u) \ge \frac{1}{2} \left[\widehat{M}(\|\nabla v\|_2^2)-\widehat{M}(\|\nabla u\|_2^2)\right]-\int_{\R^2}F(v) \mathrm{d}x, \ \
\forall \ u \in A_c, \ \forall \ v \in B_c.
\end{align}
Let $g(t)=\widehat{M}(t+s)-\widehat{M}(t)-\widehat{M}(s)$, where $t \in \R^+$.
By $(M_1)$, we get $g'(t) \ge 0$ for all $t \in \R^+$. So
\begin{align}\label{2.50}
\widehat{M}(t+s) \ge \widehat{M}(t)+\widehat{M}(s), \ \ \ \forall \ t, s \in \R^+.
\end{align}
Let
$t>1$(close to 1), $\alpha>\alpha_0$(close to $\alpha_0$) and $\tau \in (0,1)$(close to $0$) such that
$t \alpha (\|\nabla v\|_2^2+\tau\|v\|_2^2) < 4\pi$ for all $v \in B_c$.
Let $K_c<\frac{2 \pi}{\alpha_0}$ and $t'=\frac{t}{t-1}$. By Lemma \ref{lemma:20220928-hbl1}, there exists $C>0$ such that for all $v \in B_c$,
\begin{align}\label{2.52}
\int_{\R^2}|v|^q \left(e^{\alpha v^2}-1\right) \mathrm{d}x \le C \|v\|_{t'q}^q.
\end{align}
By (\ref{2-6}), (\ref{2.52}) and Lemma \ref{lemma:20220923-l1}, we get for any $\varepsilon>0$, there exists $C_\varepsilon>0$ such that
\begin{align}\label{2.53}
\int_{\R^2}F(v) \mathrm{d}x \le \varepsilon \|\nabla v\|_2^2 +C_\varepsilon\|\nabla v\|_2^{q-\frac{2}{t'}}.
\end{align}
Then by (\ref{2.49})-(\ref{2.50}), (\ref{2.53}) and $(M_1)$, we obtain that $\forall u \in A_c$, $\forall v \in B_c$,
\begin{align}\label{2.54}
I(v)-I(u) \ge \frac{1}{2}M_0 K_c -2 \varepsilon K_c- C_\varepsilon (2 K_c)^{\frac{q}{2}-\frac{1}{t'}}, \ \
\forall \ u \in A_c, \ \forall \ v \in B_c.
\end{align}
And $\forall u \in A_c$,
\begin{align}\label{212}
I(u)\ge \left(\frac{1}{2}M_0 -\varepsilon\right)\|\nabla u\|_2^2-C_\varepsilon\|\nabla u\|_2^{q-\frac{2}{t'}}.
\end{align}
Since $\frac{q}{2}-\frac{1}{t'}>1$, we can choose $\varepsilon$ and $K_c$ small such that
$\inf_{u \in B_c}I(u)>\sup_{u \in A_c}I(u)$ and $I(u)>0$ for all $u \in A_c$.
\ep

Let $I^0:=\{u \in S_{r,c}: I(u) \le 0\}$.
Define the mountain pass level
\begin{align}\label{2.55}
\gamma_c:= \inf_{h \in \Gamma_c} \max_{t \in
[0,1]}I(h (t)),
\end{align}
where
\begin{align*}
\Gamma_c:= \left\{h \in C ([0,1],S_{r,c}): h(0)\in A_c, h(1) \in I^0\right\}.
\end{align*}
Define the functional $J: H \rightarrow \R$ as follows:
\begin{align*}
J(u,s)=I(T(u,s))=\frac{1}{2}\widehat{M}(e^{2s}\|\nabla u\|_2^2)-\frac{1}{e^{Ns}}\int_{\R^N}F(e^{\frac{Ns}{2}}u)\mathrm{d}x.
\end{align*}
Define
\begin{align}\label{2.56}
\tilde{\gamma}_c:= \inf_{\tilde{h} \in \tilde{\Gamma}_c} \max_{t \in
[0,1]}J(\tilde{h} (t)),
\end{align}
where
\begin{align*}
\tilde{\Gamma}_c:= \left\{\tilde{h} \in C ([0,1],S_{r,c} \times \R): \tilde{h}(0)\in (A_c,0), \tilde{h}(1) \in (I^0,0)\right\}.
\end{align*}
Obviously, we have $\gamma_c=\tilde{\gamma}_c$. Similar to the argument of  \cite[Proposition 2.2 and Lemma 2.4]{J}, we can derive from Lemma \ref{lemma:20220928-hbl5} and Lemma \ref{lemma:20220928-hbl6} to get the following result.

\bl\lab{lemma:20220928-hbl7}There exists $\{u_n\} \subset S_{r,c}$ such that
\begin{align}\label{2.57}
I(u_n) \rightarrow \gamma_c,  \ \ I_{S_{r,c}}'(u_n) \rightarrow 0, \ \ G(u_n) \rightarrow 0.
\end{align}
\el

\subsection{Estimation of the mountain pass value}
We consider the Moser sequence of functions
\begin{align*}
\bar{\omega}_n(x)=\frac{1}{\sqrt{2\pi}}
\left\{
\begin{aligned}
 &(\log n)^\frac{1}{2},      \ \ 0\leq|x|\leq\frac{1}{n},\\
 &\frac{\log\frac{1}{|x|}}{(\log n)^\frac{1}{2}}, \ \ \frac{1}{n}\leq|x|\leq1,\\
 & 0,  \ \ \ \ \ \ \ \ \ \ \ |x|\geq 1.
\end{aligned}
\right.
\end{align*}
It is well known that
$\|\nabla \bar{\omega}_n\|_2^2=1$. By a direct calculation,
\begin{align*}
\|\bar{\omega}_n\|_2^2=&\frac{\log n}{2 n^2}+\frac{1}{ \log n}\int_{\frac{1}{n}}^1 x \log^2 x \mathrm{d}x\notag\\
=& \frac{\log n}{2 n^2}+\frac{1}{\log n}\left(\frac{1}{4}-\frac{1}{4 n^2}-\frac{\log n}{2 n^2}-\frac{\log^2 n}{2 n^2}\right)\\
=&\frac{1}{4 \log n} +o\left(\frac{1}{\log^2 n}\right).
\end{align*}
Define $\omega_n=\frac{c \bar{\omega}_n}{\|\bar{\omega}_n\|_2}$.
Then we have that
\begin{align}\label{2.3}
\|\nabla \omega_n\|_2^2=4 c^2 \log n\left(1+o\left(\frac{1}{\log n}\right)\right),
\end{align}
and
\begin{align}\label{2.4}
\omega_n(x)=\frac{\sqrt{2}c}{\sqrt{\pi}}
\left\{
\begin{aligned}
 &\log n\left(1+ o\left(\frac{1}{\log n}\right)\right),    \  \ \ 0\leq|x|\leq\frac{1}{n},\\
 &\log\frac{1}{|x|}\left(1+ o\left(\frac{1}{\log n}\right)\right), \ \ \frac{1}{n}\leq|x|\leq1,\\
 & 0,  \ \ \ \ \ \ \ \ \ \ \ \ \ \ \ \ \ \ \ \ \ \ \ \ \ \ \ \ \ \ \ \ \ |x|\geq 1.
\end{aligned}
\right.
\end{align}
Let
\begin{align}\lab{def:gn}
g_n(t)=\frac{1}{2}\widehat{M}(t^2\|\nabla \omega_n\|_2^2)-t^{-2}\int_{\R^2}F(t \omega_n)\mathrm{d}x, \ \ t \ge 0.
\end{align}
Fix $n\in \N$. By Lemma \ref{lemma:20220928-hbl5}, there exist $s_n < -1$ small and $t_n>1$ large such that $T(\omega_n,s_n) \in A_c$ and $I(T(\omega_n,t_n)) <0$. Let $h_n(t):=T(\omega_n,t t_n+ (1-t) s_n)$, where $t \in [0,1]$. Then
\begin{align}\lab{eq:20221021-e1}
\gamma_c \le \max_{t \in [0,1]}I(h_n(t)) \le \max_{t \in \R}I(T(\omega_n,t)) \le \max_{t \ge 0}g_n(t).
\end{align}

\bl\lab{lemma:20220928-hbl2} For any fixed $n\in \N$, $\max_{t \ge 0}g_n(t)>0$ is attained at some $t_n>0$.
\el

\bp By (\ref{2-6})
and Lemma \ref{lemma:20220928-hbl1}, we have that $g_n(0)=0$. Also, for $t>0$ small,
\begin{align}\label{2-7}
t^{-2}\left|\int_{\R^2}F(t \omega_n) \mathrm{d}x\right| \le & \frac{\varepsilon}{4}t^2\|\omega_n\|_4^4+
\frac{C_\varepsilon}{q}t^{q-2}\|\omega_n\|_{2q}^q\left[\int_{\R^2}(e^{2\alpha t^2 \omega_n^2}-1)\mathrm{d}x\right]^{\frac{1}{2}}\notag\\
\le & C \left(\varepsilon t^2+ C_\varepsilon t^{q-2}\right),
\end{align}
where $q>4$. By $(M_1)$, we get $\widehat{M}(s) \ge M_0 s$ for $s \in \R^+$. So
$g_n(t)>0$ for $t>0$ small.
By (\ref{2-9}) and Lemma \ref{lemma:20220929}, we have that $g_n(t)<0$ for $t>0$ large. Thus,
$\max_{t \ge 0}g_n(t)>0$ is attained at some $t_n>0$.
\ep

\bl\lab{lemma:20220928-hbl3}  For $n$ large, there holds $\displaystyle\max_{t \ge 0}g_n(t)< \frac{1}{2}\widehat{M}\left(\frac{4 \pi}{\alpha_0}\right)$.
\el

\bp By Lemma \ref{lemma:20220928-hbl2}, we get $g_n'(t_n)=0$. Let $\sigma=2\theta+4$, by $(f_4)$,
\begin{align}\label{2-10}
M(t_n^2\|\nabla \omega_n\|_2^2) t_n^2\|\nabla \omega_n\|_2^2=&-2 t_n^{-2}\int_{\R^2}F(t_n \omega_n) \mathrm{d}x+t_n^{-2}\int_{\R^2}f(t_n \omega_n)t_n \omega_n \mathrm{d}x \notag\\
 \geq & (\sigma-2) t_n^{-2} \int_{B_1(0)}F(t_n \omega_n) \mathrm{d}x.
\end{align}
By $(f_3)$, we have that
\begin{align*}
\lim_{t \rightarrow +\infty}\frac{F(t)}{t^{-2} e^{\alpha_0 t^2}}
=\lim_{t \rightarrow +\infty}\frac{f(t)}{2 \alpha_0 t^{-1} e^{\alpha_0 t^2}} \ge \frac{\beta}{2 \alpha_0}.
\end{align*}
Hence, for any $\delta>0$, there exists $t_\delta>0$ such that for $t \ge t_\delta$,
\begin{align}\label{2-11}
f(t)t \ge (\beta-\delta)e^{\alpha_0 t^2}, \ \ F(t) t^2 \ge \frac{\beta-\delta}{2 \alpha_0}e^{\alpha_0 t^2}.
\end{align}

If $\displaystyle\lim_{n \rightarrow \infty}(t_n^2 \log n) =0$, by (\ref{2.3}),
we have that $\frac{1}{2}\widehat{M}(t_n^2\|\nabla \omega_n\|_2^2)\rightarrow 0$. Then $\displaystyle\lim_{n \rightarrow \infty}g_n(t_n)=0$ and the conclusion holds. So we may assume that $\lim_{n \rightarrow \infty}(t_n^2 \log n) =l \in (0,+\infty]$. In such a case, it is trivial that
$\displaystyle\lim_{n \rightarrow \infty}(t_n \log n) =+\infty$.
By (\ref{2.3})-(\ref{2.4}) and (\ref{2-10})-(\ref{2-11}), we derive that
\begin{align*}
&M\left(4 c^2 \log n\left(1+o\left(\frac{1}{\log n}\right)\right)t_n^2\right) 4 c^2 \log n(1+o(\frac{1}{\log n}))t_n^2\\
&=M(t_n^2\|\nabla \omega_n\|_2^2) t_n^2\|\nabla \omega_n\|_2^2\\
&\geq (\sigma-2) t_n^{-2} \int_{B_{\frac{1}{n}}(0)}F(t_n \omega_n) \mathrm{d}x\\
&\ge \frac{\pi^2(\beta-\delta)(\sigma-2)}{4 \alpha_0 c^2} e^{\frac{2 \alpha_0 c^2 t_n^2}{\pi}  (\log n)^2(1+o(\frac{1}{\log n}))-2 \log n}\frac{1}{t_n^4(\log n)^2(1+o(\frac{1}{\log n}))}.
\end{align*}
If $l=+\infty$, by Lemma \ref{lemma:20220929}, we get a contradiction by the inequality above.
So $l \in (0,+\infty)$. In particular, by the inequality above again and letting $n\rightarrow +\infty$, we have that $l \in \left(0,\frac{\pi}{ \alpha_0 c^2}\right]$. If $l \in \left(0,\frac{\pi}{\alpha_0 c^2}\right)$, then
\begin{align}\label{2-12}
\lim_{n \rightarrow \infty}g_n(t_n) \le \frac{1}{2}\lim_{n \rightarrow \infty}\widehat{M}(t_n^2\|\nabla \omega_n\|_2^2)<\frac{1}{2}\widehat{M}\left(\frac{4 \pi}{\alpha_0}\right).
\end{align}

Now we consider the case $l=\frac{\pi}{\alpha_0 c^2}$.
Let
\begin{align*}
A_n:=\{x\in B_1(0):t_n\omega_n(x)\geq t_\delta\}.
\end{align*}
By (\ref{2-11}), we have that
\begin{align*}
\int_{\R^2}F(t_n \omega_n)\mathrm{d}x \geq \frac{\beta-\delta}{2 \alpha_0}\int_{A_n}t_n^{-2}\omega_n^{-2}e^{\alpha_0 t_n^2\omega_n^2}\mathrm{d}x.
\end{align*}
Let $s \in (0,\frac{1}{2})$. By \eqref{2.4}, one can see that
\beq\lab{eq:20220928-hbwhe1}
t_n\omega_n(x)\geq t_\delta, \forall \ |x| \le \frac{1}{n^s}~\hbox{uniformly for $n$ large enough}.
\eeq
So combining with \eqref{2-11}, we have that
\begin{align}\label{2-13}
t_n^{-2}\int_{\R^2}F(t_n \omega_n)\mathrm{d}x \ge \frac{\beta-\delta}{2 \alpha_0}\int_{B_{\frac{1}{n^s}}(0)}t_n^{-4}\omega_n^{-2}e^{\alpha_0 t_n^2 \omega_n^2}\mathrm{d}x.
\end{align}
By a direct calculation,
\begin{align}\label{2-14}
&\int_{B_{\frac{1}{n^s}}(0)}t_n^{-4}\omega_n^{-2}e^{\alpha_0 t_n^2 \omega_n^2}\mathrm{d}x\notag\\
&=\int_{|x| \le \frac{1}{n}}\frac{\pi e^{\frac{2\alpha_0 c^2 t_n^2 \log^2 n(1+o(\frac{1}{\log n}))}{\pi}}}{2 c^2 t_n^4 \log^2 n (1+o(\frac{1}{\log n}))}\mathrm{d}x+ \int_{\frac{1}{n} \le |x| \le \frac{1}{n^s}}\frac{\pi  e^{\frac{2 \alpha_0 c^2 t_n^2\log^2 |x|(1+o(\frac{1}{\log n}))}{\pi} }}{2 c^2 t_n^4 \log^2 |x| (1+o(\frac{1}{\log n}))}\mathrm{d}x\notag\\
&=\frac{\pi^2}{2 c^2 t_n^4}\frac{e^{\frac{2\alpha_0 c^2 t_n^2 \log^2 n(1+o(\frac{1}{\log n}))}{\pi}-2 \log n}}{\log^2 n(1+o(\frac{1}{\log n}))} +\frac{\pi^2}{ c^2 t_n^4}
\int_{\frac{1}{n}}^{\frac{1}{n^s}}\frac{x e^{\frac{2 \alpha_0 c^2 t_n^2}{\pi}\log^2 x(1+o(\frac{1}{\log n}))}}{\log^2 x(1+o(\frac{1}{\log n}))}\mathrm{d}x.
\end{align}
Let $C_n=\frac{2\alpha_0 c^2 t_n^2}{\pi}$ and $t=C_n \log \frac{1}{x}$. Then
\begin{align}\label{2-15}
&\int_{\frac{1}{n}}^{\frac{1}{n^s}}\frac{x e^{\frac{2 \alpha_0 c^2 t_n^2}{\pi}\log^2 x(1+o(\frac{1}{\log n}))}}{\log^2 x(1+o(\frac{1}{\log n}))}\mathrm{d}x\notag\\
&=C_n (1+o(\frac{1}{\log n}))\int_{s C_n \log n}^{C_n \log n}e^{-\frac{2 t}{C_n}+\frac{t^2}{C_n}(1+o(\frac{1}{\log n}))} t^{-2}\mathrm{d}t\notag\\
&\ge  \frac{1+o(\frac{1}{\log n})}{\log n}\int_s^1 e^{-2 x\log n  + C_n x^2 \log^2 n(1+o(\frac{1}{\log n})) } \mathrm{d}x.
\end{align}
Here
\begin{align*}
&\int_s^1 e^{-2 x\log n  + C_n x^2 \log^2 n(1+o(\frac{1}{\log n})) } \mathrm{d}x \notag\\
&\ge
\int_{\frac{1}{C_n \log n}}^1 e^{\left[\left(2 C_n \log^2 n (1+o(\frac{1}{\log n}))-2 \log n\right)x-C_n \log^2 n (1+o(\frac{1}{\log n}))\right]}\mathrm{d}x\notag\\
&\quad + \int_s^{\frac{1}{C_n \log n}}e^{-2x \log n }\mathrm{d}x\notag\\
&= \frac{e^{-C_n \log^2 n (1+o(\frac{1}{\log n}))}}{2 C_n \log^2 n (1+o(\frac{1}{\log n}))-2 \log n}e^{2 C_n \log^2 n (1+o(\frac{1}{\log n}))-2 \log n}\notag\\
&\quad - \frac{e^{-C_n \log^2 n (1+o(\frac{1}{\log n}))}}{2 C_n \log^2 n (1+o(\frac{1}{\log n}))-2 \log n}e^{2 \log n (1+o(\frac{1}{\log n}))-\frac{2}{C_n}}\notag\\
&\quad+\frac{1}{2 \log n}\left(e^{-2s \log n}-e^{-\frac{2}{C_n}}\right).
\end{align*}
Moreover, by $\displaystyle\lim_{n \rightarrow \infty}C_n \log n=2$, we obtain that for any $\varepsilon \in (0,1-2s)$,
there exists $N_1 \in \mathbb{N}$ such that for $n > N_1$,
\begin{align}\label{2-16}
&\int_s^1 e^{-2 x\log n  + C_n x^2 \log^2 n(1+o(\frac{1}{\log n})) } \mathrm{d}x \notag\\
&\ge \frac{e^{C_n \log^2 n (1+o(\frac{1}{\log n}))-2 \log n}}{2(1+\varepsilon)\log n}-\frac{e^{(-1+\varepsilon) \log n}}{2(1-\varepsilon)\log n}+ \frac{e^{-2s \log n} -e^{(-1+\varepsilon)\log n}}{2 \log n}\notag\\
&> \frac{e^{C_n \log^2 n (1+o(\frac{1}{\log n}))-2 \log n}}{2(1+\varepsilon)\log n}.
\end{align}
By (\ref{2-13})-(\ref{2-14}) and (\ref{2-16}), we derive that for $n$ large,
\begin{align}\label{2-17}
&t_n^{-2}\int_{\R^2}F(t_n \omega_n)\mathrm{d}x \notag\\
&\ge  \frac{(\beta-\delta)\pi^2e^{C_n \log^2 n (1+o(\frac{1}{\log n}))-2 \log n}}{4 \alpha_0 c^2 t_n^4(1+\varepsilon)\log^2 n }\notag\\
& \quad+\frac{(\beta-\delta)\pi^2(1-\varepsilon)e^{C_n \log^2 n (1+o(\frac{1}{\log n}))- 2 \log n}} {4 \alpha_0 c^2 t_n^4(1+\varepsilon)\log^2 n}\notag\\
&=\frac{(2-\varepsilon)(\beta-\delta)\pi^2 e^{C_n \log^2 n (1+o(\frac{1}{\log n}))- 2 \log n}}{4 (1+\varepsilon) \alpha_0 c^2 t_n^4 \log^2 n}.
\end{align}
Together with $(M_1)$, we obtain that for $n$ large,
\begin{align}\label{2-18}
g_n(t_n) \le & \frac{1}{2}\widehat{M}\left(4 c^2 \left(1+o\left(\frac{1}{\log n}\right)\right)t_n^2 \log n\right) \notag\\ &-\frac{(2-\varepsilon)(\beta-\delta)\pi^2 e^{C_n \log^2 n (1+o(\frac{1}{\log n}))- 2 \log n}}{4 (1+\varepsilon) \alpha_0 c^2 t_n^4 \log^2 n}\notag\\
\le & \frac{1}{2}\widehat{M}\left(4 c^2 t_n^2 \log n\right)+ o(\frac{1}{\log n}) \notag\\
&-\frac{(2-\varepsilon)(\beta-\delta)(\alpha_0 c^2-\varepsilon) e^{\frac{2 \alpha_0 c^2 t_n^2 \log^2 n (1+o(\frac{1}{\log n}))}{\pi}- 2 \log n}}{4 (1+\varepsilon)}.
\end{align}
Let
\begin{align*}
l_n(t):=\frac{1}{2}\widehat{M}\left(4 c^2 t^2\right) -\frac{(2-\varepsilon)(\beta-\delta)(\alpha_0 c^2-\varepsilon) n^{\frac{2 \alpha_0 c^2 t^2(1+o(\frac{1}{\log n})) }{\pi}- 2}}{4 (1+\varepsilon)}.
\end{align*}
So $g_n(t_n) \le \sup_{t \ge 0}l_n(t)+ o(\frac{1}{\log^2 n})$.
Moreover, there exists $t_n^*>0$ such that $\sup_{t \ge 0}l_n(t)=l_n(t_n^*)$ and $l_n'(t_n^*)=0$, from which we get
\begin{align}\label{2-19}
&M\left(4 c^2 (t_n^*)^2\right)\notag\\
&= \frac{\alpha_0(2-\varepsilon)(\beta-\delta)( \alpha_0 c^2 -\varepsilon)(1+o(\frac{1}{\log n}))}{4 \pi (1+ \varepsilon)}(\log n) n^{\frac{2 \alpha_0 c^2 (t_n^*)^2(1+o(\frac{1}{\log n})) }{\pi}- 2}.
\end{align}
By (\ref{2-18})-(\ref{2-19}), we have that
\begin{align}\label{2-20}
g_n(t_n) \leq & \frac{1}{2}\widehat{M}\left(4 c^2 (t_n^*)^2\right) -\frac{\pi}{\alpha_0 \log n (1+o(\frac{1}{\log n}))}M\left(4 c^2 (t_n^*)^2\right)+ o(\frac{1}{\log n}).
\end{align}
We claim that $\displaystyle\lim_{n \rightarrow \infty}(t_n^*)^2 = \frac{\pi}{\alpha_0 c^2}$.
If $\displaystyle\lim_{n \rightarrow \infty}(t_n^*)^2 < \frac{\pi}{\alpha_0 c^2}$, by (\ref{2-19}) and $(M_1)$, we get $M_0 \le 0$, a contradiction. If $\displaystyle\lim_{n \rightarrow \infty}(t_n^*)^2 \in \left(\frac{\pi}{\alpha_0 c^2},+\infty\right)$, by $(M_1)$, we get
\begin{align}\label{224}
\lim_{n \rightarrow +\infty}\frac{(\log n) n^{\frac{2 \alpha_0 c^2 (t_n^*)^2(1+o(\frac{1}{\log n})) }{\pi}- 2}}{M\left(4 c^2 (t_n^*)^2\right)}=+\infty,
\end{align}
a contradiction with (\ref{2-19}). Obviously, there exists $N_0 \in \mathbb{N}$ such that for $n \ge N_0$,
\begin{align*}
n^{\frac{\alpha_0 c^2 t^2}{2\pi}} \ge t^{2 \theta+2}, \ \ \ \forall \ t \ge 1.
\end{align*}
Thus, if $\displaystyle\lim_{n \rightarrow \infty}(t_n^*)^2 =+\infty$, by Lemma \ref{lemma:20220929}, we get (\ref{224}), a contradiction with (\ref{2-19}). So $\displaystyle\lim_{n \rightarrow \infty}(t_n^*)^2 =\frac{\pi}{\alpha_0 c^2}$.
Let
\begin{align*}
A_n:=\frac{4 \pi (1+ \varepsilon) M\left(4 c^2 (t_n^*)^2\right)}{\alpha_0(2-\varepsilon)(\beta-\delta)( \alpha_0 c^2 -\varepsilon)\log n(1+o(\frac{1}{\log n}))}.
\end{align*}
Then
\begin{align}\label{2-21}
(t_n^*)^2=&\frac{\pi }{\alpha_0 c^2 (1+o(\frac{1}{\log n}))}+\frac{\pi \log A_n}{2 \alpha_0 c^2 (1+o(\frac{1}{\log n}))\log n}\notag\\
\le & \frac{\pi}{\alpha_0 c^2}+o(\frac{1}{\log n})+ \frac{C \log A_n}{\log n}.
\end{align}
By (\ref{2-20}), (\ref{2-21}) and $(M_1)$, we obtain that there exists $C'>0$ such that
\begin{align}\label{2-22}
g_n(t_n) \le & \frac{1}{2}\widehat{M}\left(\frac{4 \pi}{\alpha_0}\right)
+o(\frac{1}{\log n})-\frac{C'}{\log n}.
\end{align}
By choosing $n$ large, we get $g_n(t_n) <\frac{1}{2}\widehat{M}\left(\frac{4 \pi}{\alpha_0}\right)$.
\ep

\subsection{Brezis-Lieb type results}\lab{subsection:BL-lemma}
\bl\lab{lemma:20221018-l1}
Assume that $u_n\rightharpoonup u_0$ weakly in $H^1(\R^2)$ and $\displaystyle \limsup_{n\rightarrow\infty}\|\nabla (u_n-u_0)\|_2^2<\frac{4\pi}{\alpha_0}$. Then for $\alpha>\alpha_0$ close to $\alpha_0$, we have that $\{e^{\alpha u_n^2}-1\}$ is bounded in $L^r(\R^2)$ provided $r>1$ close to $1$.
\el
\bp
We only need to prove the result for $n$ large enough.
For $\alpha>\alpha_0$ close to $\alpha_0$ and $r>1$ close to $1$, we still have that
\begin{align} \label{227}
\displaystyle \limsup_{n\rightarrow\infty}\left(r\alpha\|\nabla (u_n-u_0)\|_2^2\right)<4\pi.
\end{align}

Let $v_n=u_n-u_0$. By choosing $\sigma>0$ small enough, we have that
\begin{align*}
\int_{\R^2}(e^{\alpha u_n^2}-1)^r\mathrm{d}x
\le & \int_{\R^2}(e^{r\alpha u_n^2}-1)\mathrm{d}x\\
=&\int_{\R^2}(e^{r\alpha (u_0+v_n)^2}-1)\mathrm{d}x\\
\leq&\int_{\R^2}(e^{r\alpha [(1+\sigma)v_n^2+ (1+\frac{1}{\sigma})u_0^2]}-1)\mathrm{d}x\\
=&\int_{\R^2}\left(e^{r\alpha(1+\sigma)v_n^2}-1\right) \left(e^{r\alpha(1+\frac{1}{\sigma})u_0^2}-1\right) \mathrm{d}x +\int_{\R^2}\left(e^{r\alpha(1+\sigma)v_n^2}-1\right)\mathrm{d}x\\
&+\int_{\R^2} \left(e^{r\alpha(1+\frac{1}{\sigma})u_0^2}-1\right)\mathrm{d}x.
\end{align*}
By choosing $\eta>1$ close to $1$, we derive from \eqref{227}, Lemma \ref{lemma:20220928-hbl1} and the H\"older inequality to obtain some $C>0$ independent of $n$ such that
$$\begin{cases}
\int_{\R^2}\left(e^{r\alpha(1+\sigma)v_n^2}-1\right) \left(e^{r\alpha(1+\frac{1}{\sigma})u_0^2}-1\right) \mathrm{d}x\leq \|e^{r\alpha(1+\sigma)v_n^2}-1\|_{\eta} \|e^{r\alpha(1+\frac{1}{\sigma})u_0^2}-1\|_{\eta'}\leq C,\\
\int_{\R^2}\left(e^{r\alpha(1+\sigma)v_n^2}-1\right)\mathrm{d}x\leq C,\\
\int_{\R^2} \left(e^{r\alpha(1+\frac{1}{\sigma})u_0^2}-1\right)\mathrm{d}x\leq C,
\end{cases}$$
where $\eta'=\frac{\eta}{\eta-1}$. Hence, $\{e^{\alpha u_n^2}-1\}$ is bounded in $L^r(\R^2)$.
\ep

Then we have the following Br\'{e}zis-Lieb type result.
\bl\lab{cro:20221018-c1}
Assume that $\{u_n\}\subset H^1(\R^2)$ such that $u_n\rightharpoonup u_0$ weakly in $H^1(\R^2)$ and $\displaystyle \limsup_{n\rightarrow\infty}\|\nabla (u_n-u_0)\|_2^2<\frac{4\pi}{\alpha_0}$. Under the assumptions $(f_1)$ and $(f_2)$, we have that
\beq\lab{eq:20221018-bze1}
\int_{\R^2}F(u_n) \mathrm{d}x=\int_{\R^2}F(u_0) \mathrm{d}x+\int_{\R^2}F(u_n-u_0) \mathrm{d}x+o_n(1).
\eeq
Suppose further that $f\in C^1(\R)$, and for any $\varepsilon>0$, there exists $C_\varepsilon>0$ such that
\beq\lab{eq:20221019-sbe1}
|f'(s)|\leq \varepsilon |s|^2+C_\varepsilon |s|^{q-2} \left(e^{\alpha |s|^2}-1\right), \ \ \ \forall \ s \in \R,
\eeq
where $\alpha>\alpha_0$ and $q>4$, then we also have that
\beq\lab{eq:20221018-ze1}
\int_{\R^2}f(u_n)u_n \mathrm{d}x=\int_{\R^2}f(u_0)u_0 \mathrm{d}x+\int_{\R^2}f(u_n-u_0)(u_n-u_0) \mathrm{d}x+o_n(1).
\eeq
\el
\bp
We may assume that $u_n\rightarrow u_0$ a.e. in $\R^2$. By the mean value theorem,
there exists $\theta_n(x) \in [0,1]$ such that
\beq\lab{eq:20221019-e1}
\int_{\R^2}|F(u_n)-F(u_n-u_0)|\mathrm{d}x=\int_{\R^2}|f(u_n-(1-\theta_n)u_0) u_0| \mathrm{d}x.
\eeq
Let $\alpha>\alpha_0$ and $q>3$. By $(f_1)$-$(f_2)$, for any $\varepsilon>0$, there exists $C_\varepsilon>0$ such that
\begin{align}\label{eq:20221019-e2}
|f(u_n-(1-\theta_n)u_0) u_0|\leq &C_\varepsilon |u_n-(1-\theta_n)u_0|^{q-1} \big(e^{\alpha [u_n-(1-\theta_n)u_0]^2}-1\big) |u_0|\notag\\
&+ \varepsilon |u_n-(1-\theta_n)u_0|^2 |u_0|.
\end{align}
By Lemma \ref{lemma:20221018-l1}, $\{e^{\alpha u_n^2}-1\}$ is bounded in $L^r(\R^2)$ provided $\alpha>\alpha_0$ close to $\alpha_0$ and $r\in (1,2)$ close to $1$. Then it is easy to see that $\{e^{\alpha [u_n-(1-\theta)u_0]^2}-1\}$ is bounded in $L^{r_1}(\R^2)$ provided $1<r_1<r$.
We {\bf claim} that $\{f(u_n-(1-\theta_n)u_0)\}$ is bounded in $L^\tau(\R^2)$ for any $1<\tau<r_1$.
Indeed, by \eqref{eq:20221019-e2}, we only need to prove that $\left\{|u_n-(1-\theta_n)u_0|^{q-1} \big(e^{\alpha [u_n-(1-\theta_n)u_0]^2}-1\big)\right\}$ is bounded in $L^\tau(\R^2)$ for any $1<\tau<r_1$.
By the Sobolev embedding $H^1(\R^2)\hookrightarrow L^p(\R^2)$ for any $p\ge2$ and the H\"older inequality, we have that
\begin{align*}
&\int_{\R^2}|u_n-(1-\theta_n)u_0|^{(q-1)\tau} \big(e^{\alpha [u_n-(1-\theta_n)u_0]^2}-1\big)^\tau \mathrm{d}x\\
&\leq \left(\int_{\R^2}|u_n-(1-\theta_n)u_0|^{(q-1)\tau \frac{r_1}{r_1-\tau}}\mathrm{d}x\right)^{\frac{r_1-\tau}{r_1}} \left(\int_{\R^2} \big(e^{\alpha [u_n-(1-\theta_n)u_0]^2}-1\big)^{r_1} \mathrm{d}x\right)^{\frac{\tau}{r_1}}\\
&\leq C \|e^{\alpha [u_n-(1-\theta_n)u_0]^2}-1\|_{r_1}^{\tau}<\infty,
\end{align*}
and the claim is proved.

Noting that $u_0\in L^{\tau'}(\R^2)$ with $\tau'=\frac{\tau}{\tau-1} > 2$, we have that
$\{f(u_n-(1-\theta_n)u_0) u_0\}$ is bounded in $L^1(\R^2)$. Define $B_R^c:=\{x\in \R^2: |x|>R\}$. Then
\beq\lab{eq:20221019-e3}
\int_{B_R^c} |f(u_n-(1-\theta_n)u_0) u_0|\mathrm{d}x \leq \|f(u_n-(1-\theta_n)u_0)\|_{\tau} \|u_0\|_{L^{\tau'}(B_R^c)}\rightarrow 0,
\eeq
uniformly in $n$ as $R\rightarrow +\infty$.
Furthermore, for any $\Lambda \subset \R^2$, we have that
\beq\lab{eq:20221019-e4}
\int_{\Lambda} |f(u_n-(1-\theta_n)u_0) u_0|\mathrm{d}x \leq \|f(u_n-(1-\theta_n)u_0)\|_{L^\tau(\Lambda)} \|u_0\|_{L^{\tau'}(\Lambda)}\rightarrow 0,
\eeq
uniformly in $n$ as $\mathrm{meas}(\Lambda)\rightarrow 0$.
So $\left\{F(u_n)-F(u_n-u_0)\right\}$ possesses the uniform integrability condition, and then by the well known Vitali's convergence theorem, noting that $F(0)=0$, we obtain that
\beq
\lim_{n\rightarrow \infty}\int_{\R^2}[F(u_n)-F(u_n-u)] \mathrm{d}x=\int_{\R^2}\lim_{n\rightarrow \infty}[F(u_n)-F(u_n-u)] \mathrm{d}x=\int_{\R^2}F(u)\mathrm{d}x,
\eeq
and the result of \eqref{eq:20221018-bze1} holds.

Similarly, if $f\in C^1$ satisfies \eqref{eq:20221019-sbe1}, applying a similar argument as above, one can prove \eqref{eq:20221018-ze1}.
\ep

\bc\lab{cro:20221018-c2}
Assume that $\{u_n\}\subset H_r^1(\R^2)$ such that $u_n\rightharpoonup u_0$ weakly in $H_r^1(\R^2)$ and $\displaystyle \limsup_{n\rightarrow\infty}\|\nabla (u_n-u_0)\|_2^2<\frac{4\pi}{\alpha_0}$. Under the assumptions $(f_1)$ and $(f_2)$, we have that
\beq\lab{eq:20221018-ze2}
\int_{\R^2}f(u_n)u_n \mathrm{d}x=\int_{\R^2}f(u_0)u_0 \mathrm{d}x+o_n(1),
\eeq
and
\beq\lab{eq:20221018-bze2}
\int_{\R^2}F(u_n) \mathrm{d}x=\int_{\R^2}F(u_0)\mathrm{d}x+o_n(1).
\eeq
\ec
\bp
Recalling \eqref{2-5}, by Lemma \ref{lemma:20221018-l1}, one can prove that $\{f(u_n)\}$ is bounded in $L^\tau(\R^2)$ provided $1<\tau<\min\{r,2\}$. By the radial compact embedding, we have that $u_n\rightarrow u_0$ in $L^{\tau'}(\R^2)$, where $\tau'=\frac{\tau}{\tau-1} > 2$. Thus,
\beq\lab{eq:20221019-e5}
\int_{\R^2}|f(u_n)(u_n-u_0)|\mathrm{d}x\leq \|f(u_n)\|_{\tau} \|u_n-u_0\|_{\tau'}=o_n(1).
\eeq
On the other hand, up to a subsequence, $f(u_n)\rightharpoonup f(u_0)$ weakly in $L^\tau(\R^2)$, which implies that
\beq\lab{eq:20221019-e6}
\int_{\R^2}[f(u_n)-f(u_0)]u_0\mathrm{d}x=o_n(1).
\eeq
Hence,
\begin{align*}
\left|\int_{\R^2} f(u_n)u_n -f(u_0)u_0 \mathrm{d}x\right|\leq \left|\int_{\R^2}[f(u_n)-f(u_0)]u_0\mathrm{d}x\right|
+\left|\int_{\R^2}f(u_n)(u_n-u_0)\mathrm{d}x\right|=o_n(1)
\end{align*}
and \eqref{eq:20221018-ze2} is proved.

Basing on Lemma \ref{cro:20221018-c1}, to prove \eqref{eq:20221018-bze2}, it is sufficient to prove that
$$\int_{\R^2}F(u_n-u_0)\mathrm{d}x=o_n(1).$$
Recalling \eqref{2-6}, we only need to prove that
\beq\lab{eq:20221018-ze3}
\int_{\R^2}|u_n-u_0|^q \left(e^{\alpha (u_n-u_0)^2}-1\right)\mathrm{d}x=o_n(1).
\eeq
Indeed, by \eqref{227}, Lemma \ref{lemma:20220928-hbl1} and the H\"older inequality, there exists $C>0$ independent of $n$ such that
$$\int_{\R^2}|u_n-u_0|^q \left(e^{\alpha (u_n-u_0)^2}-1\right)\mathrm{d}x \leq \||u_n-u_0|^q\|_{r'} \|e^{\alpha (u_n-u_0)^2}-1\|_r \le C \||u_n-u_0|^q\|_{r'},$$
where $r'=\frac{r}{r-1}$. By the radial compact embedding again, we have that $u_n\rightarrow u_0$ in $L^{qr'}(\R^2)$ and thus \eqref{eq:20221018-ze3} holds.
\ep

\subsection{PSP-condition}\lab{subsection:PSP-condition}

\bl\lab{lemma:20221217} Assume that $\{u_n\} \subset H_r^1(\R^2)$ such that $u_n \rightharpoonup u_0$ weakly in $H_r^1(\R^2)$ and $\int_{\R^2}f(u_n)u_n \mathrm{d}x$ is bounded. Under the assumptions $(f_1)$-$(f_2)$ and $(f_5)$,
we have that
\begin{align}\label{2.31}
\lim_{n \rightarrow \infty}\int_{\R^2}F(u_n)\mathrm{d}x=\int_{\R^2}F(u_0)\mathrm{d}x.
\end{align}
\el

\bp By $(f_5)$, for any $\varepsilon>0$, there exists $M_\varepsilon>0$ such that
$$F(s)\leq \varepsilon f(s)s,\quad \forall |s| \geq M_\varepsilon.$$
Moreover, by $(f_1)$, we can find some $C_\varepsilon>0$ such that
\beq\lab{eq:20221213-e2}
F(s)\leq \varepsilon f(s)s +C_\varepsilon |s|^4,  \ \ \forall  s \in \R.
\eeq
Then
\begin{align*}
\int_{\R^2}F(u_n)\mathrm{d}x\leq &\varepsilon \int_{\R^2} f(u_n)u_n \mathrm{d}x +C_\varepsilon \|u_n\|_4^4.
\end{align*}
So $\{\int_{\R^2}F(u_n)\mathrm{d}x\}$ is bounded.
For any $\Omega\subset \R^2$, by $u_n\rightarrow u_0$ in $L^4(\R^2)$,
\begin{align*}
\int_{\Omega} F(u_n)\mathrm{d}x\leq & \int_{\Omega} \varepsilon  f(u_n)u_n \mathrm{d}x + C_\varepsilon \int_{\Omega} |u_n|^4 \mathrm{d}x\\
\leq &C \varepsilon + o_n(1)+C_\varepsilon \int_{\Omega} |u_0|^4 \mathrm{d}x.
\end{align*}
Then by the arbitrary of $\varepsilon$, it is easy to see that $\int_{\R^2}F(u_n)\mathrm{d}x$ satisfies the uniform absolute continuity, and thus by the Vitali's convergence theorem, we get
(\ref{2.31}).
\ep

We establish the following so called $PSP$-condition.

\bl\lab{lemma:20220928-hbl4} If $\{u_{n}\}
\subset S_{r,c}$ is a $(PSP)_{d}$ sequence with $d <\frac{1}{2}\widehat{M}\left(\frac{4 \pi}{\alpha_0}\right), d\neq 0$, i.e.,
$$I(u_n)
\rightarrow d <\frac{1}{2}\widehat{M}\left(\frac{4 \pi}{\alpha_0}\right), \ \ d\neq 0, \ \ I|'_{S_{r,c}}(u_n) \rightarrow 0~\hbox{and}~G(u_n) \rightarrow 0,$$
then $\{u_{n}\}$ converges strongly in $H_r^1(\R^2)$ up to a subsequence.
\el

\bp Define the functional $\varphi(u)=\frac{1}{2}\|u\|_2^2$, where $u \in H_r^1(\R^2)$.
By $I|'_{S_{r,c}}(u_n) \rightarrow 0$, we get there exists $\{\lambda_n\} \subset \R$ such that
\begin{align}\label{2.33}
\|I'(u_n) - \lambda_n \varphi'(u_n)\| \rightarrow 0.
\end{align}
By $I(u_n) \rightarrow d$, $G(u_n) \rightarrow 0$ and $(f_4)$,
\begin{align*}
d \geq \frac{1}{2}\widehat{M}(\|\nabla u_n\|_2^2)-\frac{1}{2(\theta+1)}M(\|\nabla u_n\|_2^2)\|\nabla u_n\|_2^2+o_n(1).
\end{align*}
So by $(M_2)$, we get $\|\nabla u_n\|_2$ is bounded. Moreover, $\|u_n\|_{H^1}$, $\int_{\R^2}F(u_n)\mathrm{d}x$,
$\int_{\R^2}f(u_n)u_n\mathrm{d}x$ and $|\lambda_n|$ are bounded. Going to a subsequence, we may assume that $u_n \rightharpoonup u_0$ weakly in $H_r^1(\R^2)$ and $\lambda_n \rightarrow \lambda_0$ as $n \rightarrow \infty$.

We consider two cases.

{\bf Case 1:} $u_n \rightharpoonup 0$ weakly in $H_r^1(\R^2)$.

By Lemma \ref{lemma:20221217}, we have that $\int_{\R^2}F(u_n)\mathrm{d}x \rightarrow 0$. Then
$d =  \frac{1}{2}\widehat{M}\left(\lim_{n \rightarrow \infty}\|\nabla u_n\|_2^2\right)$.
Together with $(M_1)$, we get
\begin{align}\label{2.25}
\lim_{n \rightarrow \infty}\|\nabla u_n\|_2^2 <\frac{4 \pi}{\alpha_0}.
\end{align}
By Corollary \ref{cro:20221018-c2}, we have that
$\lim_{n\rightarrow\infty}\int_{\R^2}f(u_n)u_n\mathrm{d}x=0$.
Moreover, by $G(u_n)\rightarrow 0$, we obtain that $\|\nabla u_n\|_2 \rightarrow 0$, a contradiction with $d\neq 0$.

{\bf Case 2:} $u_n \rightharpoonup u_0 \ne 0$ weakly in $H_r^1(\R^2)$.

Define
\begin{align}\lab{eq:20221018-e1}
D=\lim_{n \rightarrow \infty}\frac{\widehat{M}(\|\nabla u_n\|_2^2)}{\|\nabla u_n\|_2^2}, \ \ E=\lim_{n \rightarrow \infty}M(\|\nabla u_n\|_2^2).
\end{align}
Let $v_n=u_n-u_0$.
By Lemma \ref{lemma:20221217}, we have that
\begin{align}\label{2.32}
d \geq   \frac{D}{2}\left(\lim_{n \rightarrow \infty}\|\nabla v_n\|_2^2+\|\nabla u_0\|_2^2\right)-\int_{\R^2}F(u_0)\mathrm{d}x.
\end{align}
By (\ref{2.33}), $G(u_n) \rightarrow 0$ and Lemma \ref{lemma:20221217}, we derive that
in the weak sense,
\begin{align}\label{2.35}
-E \Delta u_0-\lambda_0 u_0= f(u_0)~\hbox{in}~\R^2,
\end{align}
with
\begin{align}\label{2.34}
\lambda_0=-\frac{2}{c^2}\int_{\R^2}F(u_0)\mathrm{d}x<0.
\end{align}
Then we have the Pohozaev identity: $\lambda_0 \|u_0\|_2^2=-2\int_{\R^2}F(u_0)\mathrm{d}x$. So $u_0\in S_{r,c}$ and thus $u_n\rightarrow u_0$ in $L^2(\R^2)$.
Moreover,
\begin{align}\label{2.36}
E \|\nabla u_0\|_2^2 + 2 \int_{\R^2}F(u_0)\mathrm{d}x-\int_{\R^2}f(u_0)u_0\mathrm{d}x=0.
\end{align}
By (\ref{2.36}), $(M_2)$ and $(f_4)$,
\begin{align}\label{2.37}
\frac{D}{2}\|\nabla u_0\|_2^2-\int_{\R^2}F(u_0)\mathrm{d}x=&\frac{D}{E}\int_{\R^2}
\left[\frac{1}{2}f(u_0)u_0-\left(1+\frac{E}{D}\right)F(u_0)\right] \mathrm{d}x\notag\\
\ge &\frac{D}{2E}\int_{\R^2}\left[f(u_0)u_0-(2\theta+4)F(u_0)\right]\mathrm{d}x \ge 0.
\end{align}
By $(M_1)$, we get $\frac{\widehat{M}(t)}{t}$ is increasing for $t>0$.
So by (\ref{2.32}) and (\ref{2.37}),
\begin{align}\label{2.38}
d \geq   \frac{D}{2}\lim_{n \rightarrow \infty}\|\nabla v_n\|_2^2
\ge \lim_{n \rightarrow \infty}\widehat{M}(\|\nabla v_n\|_2^2).
\end{align}
Moreover,
\begin{align}\label{2.39}
\lim_{n \rightarrow \infty}\|\nabla v_n\|_2^2 <\frac{4 \pi}{\alpha_0}.
\end{align}
Then by Corollary \ref{cro:20221018-c2}, we have that
\beq\lab{eq:20221018-e3}
\lim_{n \rightarrow \infty}\int_{\R^2}f(u_n)u_n\mathrm{d}x= \int_{\R^2}f(u_0)u_0\mathrm{d}x.
\eeq
Now, by \eqref{2.31}, \eqref{2.36}, \eqref{eq:20221018-e3} and the fact of $G(u_n)\rightarrow 0$, we have that $\|\nabla u_n\|_2^2\rightarrow \|\nabla u_0\|_2^2$, which implies that $u_n\rightarrow u_0$ in $D^{1,2}(\R^2)$.
Hence, $u_n\rightarrow u_0$ in $H_r^1(\R^2).$
 \ep

\subsection{Proof of  Theorem  \ref{Theorem 1.1}}\lab{subsection:proof-th1-1}
By Lemma \ref{lemma:20220928-hbl3} and \eqref{eq:20221021-e1}, we have that $\gamma_c< \frac{1}{2}\widehat{M}\left(\frac{4 \pi}{\alpha_0}\right)$. By Lemma \ref{lemma:20220928-hbl7}, there exists a $(PSP)_{\gamma_c}$-sequence. By Lemma \ref{lemma:20220928-hbl6}, we have that $\gamma_c>0$. Then by Lemma \ref{lemma:20220928-hbl4}, we get the result.
\hfill $\square$

\br\lab{remark:20221018-xr1}
For any $0<a\leq b<+\infty$, define
\beq\lab{eq:20221018-xe1}
\mathcal{U}_{r,a}^{b}:=\left\{u\in S_{r,c}: I(u)=\gamma_c, I'(u)=0, c \in [a,b]\right\}.
\eeq
which is compact in $H_r^1(\R^2)$.
Indeed, it is standard to prove that $\gamma_c$ depends continuously on $c>0$. On the other hand, by Lemma \ref{lemma:20220928-hbl3}, $\gamma_c< \frac{1}{2}\widehat{M}\left(\frac{4 \pi}{\alpha_0}\right)$ for all $c\in [a,b]$.
Then  similar to the proof of Lemma \ref{lemma:20220928-hbl4}, one can prove the conclusion.
\er

\s{The case $N \ge 5$}\lab{sec:caseNgeq5}

\renewcommand{\theequation}{3.\arabic{equation}}
\subsection{Some preliminaries}\lab{subsection:N-5}

Define the best Sobolev constant:
\begin{align}\label{4.1}
S:=\inf_{u \in D^{1,2}(\R^N) \setminus \{0\}}\frac{\int_{\R^N}|\nabla u|^2 \mathrm{d}x}{\left(\int_{\R^N}|u|^{2^*}\mathrm{d}x\right)^{\frac{2}{2^*}}}.
\end{align}
Consider the case of $M(t)=a+bt$, $a$, $b>0$ and $f(u)=|u|^{p-2}u$, $p \in (2,2^*)$. In such a case, the corresponding energy functional given by \eqref{eq:def-energy-functional} is
\begin{align*}
I(u)=\frac{a}{2}\|\nabla u\|_2^2+\frac{b}{4}\|\nabla u\|_2^4-\frac{1}{p}\|u\|_p^p-\frac{1}{2^*}\|u\|_{2^*}^{2^*}.
\end{align*}
By (\ref{4.1}) and Lemma \ref{lemma:20220923-l1}, we have that
\begin{align}\label{4.2}
I(u) \ge \frac{a}{2}\|\nabla u\|_2^2 + \frac{b}{4} \|\nabla u\|_2^4-\frac{c^{p-\frac{N(p-2)}{2}}}{2 \|Q\|_2^{p-2}}\|\nabla u\|_2^{\frac{N(p-2)}{2}}-\frac{1}{2^* S^{\frac{2^*}{2}}}\|\nabla u\|_2^{2^*}
\end{align}
holds for all $u\in S_{c}$.

\bl\lab{lemma:20220923-hhl1}
Let $\kappa_1$, $\kappa_2>0$, $0<p_1<p_2<+\infty$, $\kappa_3\geq 0$, $\kappa_4\geq 0$ and $(\kappa_3,\kappa_4)\neq (0,0)$. For $p_3$, $p_4\in (p_1,p_2)$, we define
\beq\lab{eq:20220924-e1}
\Omega_{\kappa_1,\kappa_2,\kappa_3,\kappa_4}^{p_1,p_2,p_3,p_4}:=\inf_{t>0}\frac{\kappa_1 t^{p_1}+\kappa_2 t^{p_2}}{\kappa_3 t^{p_3}+\kappa_4 t^{p_4}}.
\eeq
Then $\Omega_{\kappa_1,\kappa_2,\kappa_3,\kappa_4}^{p_1,p_2,p_3,p_4}>0$ is attained.
Furthermore,
\begin{itemize}
\item[(i)]$\Omega_{\kappa_1,\kappa_2,\kappa_3,\kappa_4}^{p_1,p_2,p_3,p_4}$ is continuous with respect to $\kappa_i\in (0,+\infty)$, $i=1,2,3,4.$
\item[(ii)]$\Omega_{\kappa_1,\kappa_2,\kappa_3,\kappa_4}^{p_1,p_2,p_3,p_4}$ increases strictly for $\kappa_i\in (0,+\infty)$, $i=1,2$.
\item[(iii)]$\Omega_{\kappa_1,\kappa_2,\kappa_3,\kappa_4}^{p_1,p_2,p_3,p_4}$ decreases strictly for $\kappa_i\in (0,+\infty)$, $i=3,4$. In particular,
    \beq\lab{eq:20220924-zbe2}
    \lim_{\kappa_i\rightarrow +\infty}\Omega_{\kappa_1,\kappa_2,\kappa_3,\kappa_4}^{p_1,p_2,p_3,p_4}=0, i=3,4.
    \eeq
\end{itemize}
\el
\bp
By a direct calculation, we get the result.
\ep

\br\lab{remark:20220923-wr1}
By the definition of $\Omega_{\kappa_1,\kappa_2,\kappa_3,\kappa_4}^{p_1,p_2,p_3,p_4}$, one can see that
$$\kappa_1 t^{p_1}+\kappa_2 t^{p_2}\geq \Omega_{\kappa_1,\kappa_2,\kappa_3,\kappa_4}^{p_1,p_2,p_3,p_4}\big(\kappa_3 t^{p_3}+\kappa_4 t^{p_4}\big), \forall \ t  \geq 0.$$
In particular, for $\kappa_1=A>0$, $\kappa_2=B>0$, $\kappa_3=0$, $\kappa_4=1$ and $p_1=2$, $p_2=4$, $p_4=q\in (2,4)$, a direct computation shows that
\beq\lab{eq:20220926-e1}
\Omega_{A,B,0,1}^{2,4,p_3,q}=2 (q-2)^{-(\frac{q}{2}-1)}(4-q)^{-(2-\frac{q}{2})} A^{2-\frac{q}{2}} B^{\frac{q}{2}-1}.
\eeq
And we note that
\beq\lab{eq:20220926-e2}
\Omega_{A,B,0,\kappa_4}^{2,4,p_3,q}=\frac{1}{\kappa_4}\Omega_{A,B,0,1}^{2,4,p_3,q}
=\frac{1}{\kappa_4}\left[2 (q-2)^{-(\frac{q}{2}-1)}(4-q)^{-(2-\frac{q}{2})} A^{2-\frac{q}{2}} B^{\frac{q}{2}-1}\right].
\eeq
In the present paper,
we take $(\kappa_1,\kappa_2,\kappa_3,\kappa_4)=(\frac{a}{2}, \frac{b}{4}, 0, \frac{1}{2^* S^{\frac{2^*}{2}}})$ or $(\kappa_1,\kappa_2,\kappa_3,\kappa_4)=(a, b, 0,S^{-\frac{2^*}{2}} )$ as applications. We also note that
\beq\lab{eq:20220926-e3}
\Omega_{\frac{a}{2},\frac{b}{4},0,\frac{1}{2^* S^{\frac{2^*}{2}}}}^{2,4,p_3,2^*}
=\left(2^{\frac{2^*}{2}-1} 4^{\frac{2-2^*}{2}} 2^*\right) \Omega_{a,b,0,\frac{1}{ S^{\frac{2^*}{2}}}}^{2,4,p_3,2^*}>\Omega_{a,b,0,\frac{1}{ S^{\frac{2^*}{2}}}}^{2,4,p_3,2^*}
\eeq
since
$$2^{\frac{2^*}{2}-1} 4^{\frac{2-2^*}{2}} 2^*=\frac{2^*}{2^{\frac{2^*-2}{2}}}>1~\hbox{for}~N\geq 5.$$
In particular, by \eqref{eq:20220926-e2},
$\displaystyle\Omega_{a,b,0,S^{-\frac{2^*}{2}}}^{2,4,p_3,2^*}>1$
equivalents to
\beq\lab{eq:20220924-ze1}
\left(\frac{2a}{4-2^*}\right)^{\frac{4-2^*}{2}} \left(\frac{2b}{2^*-2}\right)^{\frac{2^*-2}{2}}>\frac{1}{S^{\frac{2^*}{2}}}.
\eeq
\er\hfill$\Box$

\bl\lab{lemma:property-Ic}
For $c>0$, let $I_c$ be defined by \eqref{def:Ic}.
If $N\geq 5$ and $p\in (2,2^*)$, then
\begin{itemize}
\item[(i)]$-\infty<I_c\leq 0$ for all $c>0$;
\item[(ii)] $I_c$ is non-increasing on $c \in (0,+\infty)$;
\item[(iii)]$I_c$ is continuous on $c \in (0,+\infty)$.
\end{itemize}
\el
\bp
(i) Since $2^*<2+\frac{8}{N}$ for $N\geq 5$, we have that $\frac{N(p-2)}{2}<4$. So by \eqref{4.2}, we get
\beq\lab{eq:20220923-zbe1}
\lim_{\|\nabla u\|_2 \rightarrow +\infty}I(u)=+\infty.
\eeq
Moreover, $I_c> -\infty$ for all $c>0$. In particular, fix $u \in S_c$ and let $u^t(x):=t^{\frac{N}{2}} u(tx)$, $t>0$.
Then
\begin{align*}
I_c \le I(u^t) \le \frac{a}{2} t^2 \|\nabla u\|_2^2 + \frac{b}{4} t^4 \|\nabla u\|_2^4.
\end{align*}
Let $t \rightarrow 0+$, we get $I_c \le 0$.

(ii) Let $c_2>c_1$. By the definition of $I_{c_1}$, there exists $\{u_n\}\subset S_{c_1}$ such that
$$I_{c_1}\leq I(u_n)\leq I_{c_1}+\frac{1}{n}.$$
Let $$v_n(x):=\left(\frac{c_1}{c_2}\right)^{\frac{N-2}{2}}u_n\left(\frac{c_1}{c_2}x\right).$$
Then $v_n\in S_{c_2}$. In particular, it holds that
$$\|\nabla v_n\|_2=\|\nabla u_n\|_2, \ \|v_n\|_{2^*}=\|u_n\|_{2^*}~\hbox{and}~\|v_n\|_p^p=\left(\frac{c_2}{c_1}\right)^{N-\frac{N-2}{2}p}\|u_n\|_p^p.$$
So
\begin{align}\label{eq:20220922-e1}
I_{c_2}\leq I(v_n)=I(u_n)+\frac{1}{p}\left[1-\left(\frac{c_2}{c_1}\right)^{N-\frac{N-2}{2}p}\right]\|u_n\|_p^p<I(u_n)\leq I_{c_1}+\frac{1}{n}.
\end{align}
By letting $n\rightarrow +\infty$, we obtain that $I_{c_2}\leq I_{c_1}$.

(iii) For any $c>0$ and any sequence $\{c_n\}\subset \R^+$ with $c_n \rightarrow c$ as $n \rightarrow \infty$, by the definition of $I_{c_n}$, there exists $\{u_n\} \subset S_{c_n}$ such that $I_{c_n} \le I(u_n) \le I_{c_n}+\frac{1}{n}$. Let $v_n:=u_n(\theta_n^{-\frac{2}{N}}x)$, where $\theta_n=\frac{c^2}{c_n^2}=1+o_n(1)$.
Then $v_n \in S_{c}$. Moreover,
\begin{align*}
\|\nabla v_n\|_2^2=\theta_n^{2-\frac{4}{N}}\|\nabla u_n\|_2^2, \ \ \|v_n\|_p^p=\theta_n^2 \|u_n\|_p^p~\hbox{and}~
\|v_n\|_{2^*}^{2^*}=\theta_n^2 \|u_n\|_{2^*}^{2^*}.
\end{align*}
By $I_c\leq 0$ for all $c>0$ and \eqref{eq:20220923-zbe1}, one can see that $\{u_n\}$ is bounded in $H^1(\R^N)$.
So
\begin{align*}
I_{c} \le I(v_n)=\frac{a}{2}\theta_n^{2-\frac{4}{N}}\|\nabla u_n\|_2^2+\frac{b}{4}\theta_n^{4-\frac{8}{N}}\|\nabla u_n\|_2^4-\frac{1}{p}\theta_n^2 \|u_n\|_p^p-\frac{1}{2^*}\theta_n^2 \|u_n\|_{2^*}^{2^*}
=I(u_n)+o_n(1).
\end{align*}
Let $n \rightarrow \infty$,
we get $\displaystyle I_c \le \liminf_{n \rightarrow \infty}I_{c_n}$.

Similarly, for any $n\in \N$, there exists some $u_n\in S_c$ such that
$\displaystyle I_c\leq I(u_n)\leq I_c+\frac{1}{n}$.
Then $v_n(x):=u_n(\theta_n^{-\frac{2}{N}}x)\in S_{c_n}$, where $\theta_n=\frac{c_n^2}{c^2}=1+o_n(1)$. So
$I_{c_n}\leq I(v_n)=I(u_n)+o_n(1)$, and thus $\displaystyle\limsup_{n \rightarrow \infty}I_{c_n} \le I_c$.
\ep

\bc\label{cro:20220922-c1}
Let $N\geq 5$ and $p \in \left(2,2^*\right)$.
If there exists $0<\underbar{c}<\bar{c}<+\infty$ such that $I_{\underbar{c}}=I_{\bar{c}}$, then $I_c$ is not attained provided $c\in [\underbar{c}, \bar{c})$.
\ec
\bp
By (ii) of Lemma \ref{lemma:property-Ic}, we have that $I_c\equiv I_{\bar{c}}, \forall c\in [\underbar{c},\bar{c}]$. Suppose that there exists some $c\in [\underbar{c},\bar{c})$ which is attained by some $u\in S_c$. That is, $I(u)=I_c=I_{\bar{c}}$. Define
$$v(x):=\left(\frac{c}{\bar{c}}\right)^{\frac{N-2}{2}}u\left(\frac{c}{\bar{c}}x\right),$$
then $v\in S_{\bar{c}}$ and thus
$$I_{\bar{c}}\leq I(v)=I(u)+\frac{1}{p}\left[1-\left(\frac{\bar{c}}{c}\right)^{N-\frac{N-2}{2}p}\right]\|u\|_p^p<I(u)=I_{c},$$
a contradiction.
\ep

\bl\lab{lemma:20220927-zbbl1}
Let $N\geq 5$ and $p \in \left(2,2^*\right)$. Then $I_c<0$ for $c>0$ large enough. And thus we can define that
\begin{align}\lab{def:c1}
c_1:=\inf\{c>0: I_c<0\}\in [0,+\infty).
\end{align}
Furthermore, if $c_1>0$, then $I_c\equiv 0$ for $c\in (0,c_1]$ and $I_c$ is not attained provided $c\in (0,c_1)$.
\el

\bp
Recalling Lemma \ref{lemma:20220923-l1},
\begin{align}\label{4.8}
\|\nabla Q\|_2^2=\|Q\|_2^2=\frac{2}{p}\|Q\|_p^p.
\end{align}
Let $Q_t(x):=\frac{c t^{\frac{N}{2}} Q(tx)}{\|Q\|_2}$, then by \eqref{4.8},
\begin{align}\label{4.9}
I(Q_t)=\frac{a}{2} c^2 t^2 + \frac{b}{4} c^4 t^4-\frac{c^p t^{\frac{Np}{2}-N}}{2 \|Q\|_2^{p-2}}-\frac{c^{2^*} \|Q\|_{2^*}^{2^*} t^{2^*}}{2^* \|Q\|_2^{2^*}}.
\end{align}

If $2<p<2+\frac{4}{N}$, then $0<\frac{Np}{2}-N<2$. So for any $c>0$, it is easy to see that $I(Q_t)<0$ for $t>0$ small enough,  and thus
$I_c<0$ for any $c>0$. Hence, $c_1=0$.

If $p=2+\frac{4}{N}$, put $s=ct$, then
$$I(Q_t)=\left[\frac{a}{2}-\frac{c^{\frac{4}{N}}}{2\|Q\|_{2}^{p-2}}\right]s^2+\frac{b}{4}s^4-
\frac{\|Q\|_{2^*}^{2^*} }{2^* \|Q\|_2^{2^*}} s^{2^*}.$$
Then $\frac{a}{2}-\frac{c^{\frac{4}{N}}}{2\|Q\|_{2}^{p-2}}<0$ provided $c>c_0:=\left(a \|Q\|_2^{p-2}\right)^{\frac{N}{4}}$. So for $t>0$ small enough, we have that $s>0$ small enough and then $I_c\leq I(Q_t)<0$ for $c>c_0$. Hence, $c_1\in [0,+\infty)$ is also well defined.

If $p\in (2+\frac{4}{N}, 2^*)$, we rewrite
\beq\lab{eq:20220924-zbe1}
I(Q_t)=\frac{a}{2}s^2+\frac{b}{4}s^4-\frac{t^{-(N+p-\frac{Np}{2})}}{2\|Q\|_{2}^{p-2}} s^p -
\frac{\|Q\|_{2^*}^{2^*} }{2^* \|Q\|_2^{2^*}} s^{2^*}, \hbox{with}~s=ct.
\eeq
Noting that in such a case, Lemma \ref{lemma:20220923-hhl1} is applied and by \eqref{eq:20220924-zbe2}, there exists $\kappa_3>0$ large enough such that
\beq\lab{eq:20220924-zbe3}
\Omega_{\frac{a}{2},\frac{b}{4},\kappa_3,\frac{\|Q\|_{2^*}^{2^*} }{2^* \|Q\|_2^{2^*}}}^{2,4,p,2^*}<1.
\eeq
In particular, by the definition of $\Omega_{\frac{a}{2},\frac{b}{4},\kappa_3,\frac{\|Q\|_{2^*}^{2^*} }{2^* \|Q\|_2^{2^*}}}^{2,4,p,2^*}$, there exists some $s_0>0$ such that
\beq\lab{eq:20220924-zbe4}
\frac{a}{2}s_0^2+\frac{b}{4}s_0^4-\kappa_3 s_0^p -
\frac{\|Q\|_{2^*}^{2^*} }{2^* \|Q\|_2^{2^*}} s_0^{2^*}<0.
\eeq
Since $p<2^*$, there exists some $c_0>1$ large enough such that
$\displaystyle\frac{(\frac{s_0}{c_0})^{-(N+p-\frac{Np}{2})}}{2\|Q\|_{2}^{p-2}}>\kappa_3$.
Then for any $c>c_0$, by $t_0=\frac{s_0}{c}<\frac{s_0}{c_0}$, we have that
$\displaystyle\frac{t_0^{-(N+p-\frac{Np}{2})}}{2\|Q\|_{2}^{p-2}}>\kappa_3,$
and thus
$$I_c\leq I(Q_{t_0})<\frac{a}{2}s_0^2+\frac{b}{4}s_0^4-\kappa_3 s_0^p -
\frac{\|Q\|_{2^*}^{2^*} }{2^* \|Q\|_2^{2^*}} s_0^{2^*}<0, \ \forall  c>c_0.$$

It is trivial that $c_1\geq 0$. If $c_1>0$, by Lemma \ref{lemma:property-Ic}, one can see that
$$I_c\begin{cases}
\equiv 0,\quad &\hbox{if}~c\in (0,c_1],\\
<0,\quad&\hbox{if}~c>c_1.
\end{cases}$$
And thus, by Corollary \ref{cro:20220922-c1}, $I_c$ is not attained for $c\in (0,c_1)$.
\ep

To establish the existence result, we need the following $PSP$-condition.
\bl\lab{lemma:20220923-zl1} Let $N\geq 5$, $p \in \left(2,2^*\right)$ and $\displaystyle\Omega_{a,b,0,\frac{1}{S^{\frac{2^*}{2}}}}^{2,4,p,2^*}>1$.
Assume that $\{u_n\}\subset S_{r,c}$ is a $(PSP)_{d}$ sequence with $d\neq 0$, i.e.,
$$I(u_n)
\rightarrow d \ne 0, \ \ I|'_{S_{r,c}}(u_n) \rightarrow 0~\hbox{and}~ G(u_n) \rightarrow 0,$$
then $\{u_{n}\}$ converges strongly in $H_r^1(\R^N)$ up to a subsequence.
\el

\bp
By \eqref{eq:20220923-zbe1}, we get $\|\nabla u_n\|_2$ is bounded. And thus, $\{u_n\}$ is bounded in $H_r^1(\R^N)$. Assume that $u_n \rightharpoonup u_0$ weakly in $H_r^1(\R^N)$. We assume that $\displaystyle E=\lim_{n \rightarrow \infty}\|
\nabla u_n\|_2^2$. Noting that there exists $\{\lambda_n\}\subset \R$ such that
\beq
I'(u_n)=I|'_{S_{r,c}}(u_n)+\lambda_n u_n.
\eeq
So by $I|'_{S_{r,c}}(u_n) \rightarrow 0$, we get
 \beq
 I'(u_n)u_n=o_n(1)\|u_n\|_{H^1}+\lambda_n\|u_n\|_2^2.
 \eeq
 Furthermore, by the boundedness of $\{u_n\}$ in $H_r^1(\R^N)$ and $\{u_n\}\subset S_{r,c}$, we obtain that $|\lambda_n|$ is bounded. So we can assume that  $\lambda_n \rightarrow \lambda_0$. Now, by the facts of
\beq\label{eq:20220923-xbe1}
\begin{cases}
&I'(u_n) - \lambda_n u_n \rightarrow 0,\\
&\lambda_n \rightarrow \lambda_0,\\
&\|u_n\|_2^2=c, \|\nabla u_n\|_2^2\rightarrow E,\\
&u_n\rightharpoonup u_0~\hbox{weakly in $H_r^1(\R^N)$ and $u_n\rightarrow u_0$ in $L^p(\R^N)$},
\end{cases}
\eeq
it is easy to check that $u_0\in H_r^1(\R^N)$ weakly solves
\beq\lab{eq:20220923-zbe2}
-(a+bE)\Delta u_0-\lambda_0 u_0=|u_0|^{p-2}u_0+|u_0|^{2^*-2}u_0~\hbox{in}~\R^N.
\eeq
Testing \eqref{eq:20220923-zbe2} by $u_0$, we have that
\beq\lab{eq:20220923-zbe3}
(a+bE)\|\nabla u_0\|_2^2-\lambda_0\|u_0\|_2^2=\|u_0\|_p^p+\|u_0\|_{2^*}^{2^*}.
\eeq
Since $u_0$ satisfies the following Pohozaev identity:
\beq\lab{eq:20220923-zbe4}
(N-2)(a+bE)\|\nabla u_0\|_2^2-N\lambda_0\|u_0\|_2^2=2N\left[\frac{1}{p}\|u_0\|_p^p+\frac{1}{2^*}\|u_0\|_{2^*}^{2^*}\right],
\eeq
we obtain that
\beq\lab{eq:20220923-zbe5}
\lambda_0\|u_0\|_2^2=\left[\frac{N}{2}-\frac{N}{p}-1\right]\|u_0\|_p^p.
\eeq
On the other hand, by $I|'_{S_{r,c}}(u_n) \rightarrow 0$ and $G(u_n) \rightarrow 0$, we have that
\beq\lab{eq:20220923-zbe6}
\begin{cases}
(a+b\|\nabla u_n\|_2^2)\|\nabla u_n\|_2^2=\lambda_n\|u_n\|_2^2+\|u_n\|_p^p+\|u_n\|_{2^*}^{2^*}+o_n(1),\\
(a+b\|\nabla u_n\|_2^2)\|\nabla u_n\|_2^2=\frac{(p-2)N}{2p}\|u_n\|_p^p+\|u_n\|_{2^*}^{2^*}+o_n(1).
\end{cases}
\eeq
So by the boundedness of $\{u_n\}$ in $H_r^1(\R^N)$ again, we obtain that
\beq\lab{eq:20220923-zbe7}
\lambda_n c^2=\lambda_n \|u_n\|_2^2=\left[\frac{N}{2}-\frac{N}{p}-1\right]\|u_n\|_p^p+o_n(1).
\eeq
By $u_n\rightarrow u_0$ in $L^p(\R^N)$, \eqref{eq:20220923-zbe5} and \eqref{eq:20220923-zbe7}, we get $\lambda_n c^2\rightarrow \lambda_0\|u_0\|_2^2$. Furthermore, by $\lambda_n\rightarrow \lambda_0$, we obtain that
\beq\lab{eq:20220923-zbe8}
\lambda_0(\|u_0\|_2^2-c^2)=0.
\eeq

If $u_0=0$, we have that $\lambda_0=0$ and $\|u_n\|_p^p\rightarrow 0$. So by \eqref{eq:20220923-zbe6} and the Sobolev inequality,
\begin{align}\lab{eq:20220923-zbe10}
\left(a+ b \|\nabla u_n\|_2^2\right)\|\nabla u_n\|_2^2 =\|u_n\|_{2^*}^{2^*}+ o_n(1) \le \frac{\|\nabla u_n\|_{2}^{2^*}}{S^{\frac{2^*}{2}}}+ o_n(1).
\end{align}
By letting $n\rightarrow +\infty$ in \eqref{eq:20220923-zbe10}, we obtain that
$$a (\sqrt{E})^2 +b(\sqrt{E})^4\leq S^{-\frac{2^*}{2}} (\sqrt{E})^{2^*}.$$
So by $\displaystyle\Omega_{a,b,0,S^{-\frac{2^*}{2}}}^{2,4,p,2^*}>1$, we obtain that
 $\|\nabla u_n\|_2^2 \rightarrow E=0$. By the Sobolev inequality again, we also have that $\|u_n\|_{2^*}^{2^*}\rightarrow 0$.
Hence, $\displaystyle d=\lim_{n\rightarrow \infty}I(u_n)=0$, a contradiction.

Now, by $u_0\neq 0$, \eqref{eq:20220923-zbe7} implies that
$$\lambda_0=\frac{1}{c^2}\left[\frac{N}{2}-\frac{N}{p}-1\right]\|u_0\|_p^p.$$
In particular, by $p<2^*$, we have that $\lambda_0<0$. And by \eqref{eq:20220923-zbe8}, we obtain that $u_0\in S_c$ and thus $u_n\rightarrow u_0$ in $L^2(\R^N)$. Let $v_n:=u_n-u_0$.
By \eqref{eq:20220923-zbe3}, \eqref{eq:20220923-zbe6} and the Brezis-Lieb lemma in \cite{Willem}, we have that
\begin{align}\label{4.17}
\left(a+ b \|\nabla v_n\|_2^2\right)\|\nabla v_n\|_2^2\leq &\left(a+ b E\right)\|\nabla v_n\|_2^2 +o_n(1) \notag\\
 =&\|v_n\|_{2^*}^{2^*}+ o_n(1)\leq \frac{\|\nabla v_n\|_{2}^{2^*}}{S^{\frac{2^*}{2}}}+o_n(1).
\end{align}
By $\displaystyle\Omega_{a,b,0,S^{-\frac{2^*}{2}}}^{2,4,p,2^*}>1$ again, we obtain that $\|\nabla v_n\|_2\rightarrow 0$. That is, $u_n\rightarrow u_0$ in $D^{1,2}(\R^N)$. Hence, we obtain that $u_n \rightarrow u_0$ in $H_r^1(\R^N)$.
\ep

We can also have the following non-existence result.

\bt\lab{th:20220928-t1}  Let $N\geq 5$, $\displaystyle\Omega_{a,b,0,\frac{1}{S^{\frac{2^*}{2}}}}^{2,4,p,2^*}>1$ and $p \in \left[2+\frac{4}{N},2^*\right)$. Then there exists $c_0>0$ such that problem \eqref{1.5}-\eqref{1.6} has no normalized solution provided $c < c_0$.
\et

\bp If $u\in H^1(\R^N)$ solves \eqref{1.5}-\eqref{1.6}, it holds the following Pohozaev identity:
\begin{align}\label{4.31}
2a\|\nabla u\|_2^2+2b\|\nabla u\|_2^4-2 \|u\|_{2^*}^{2^*}= \left(N-\frac{2N}{p}\right)\|u\|_p^p.
\end{align}
Moreover, by (\ref{4.1}) and Lemma \ref{lemma:20220923-l1},
\begin{align}\label{4.32}
2a\|\nabla u\|_2^2+2b\|\nabla u\|_2^4 \le \frac{N(p-2) c^{p-\frac{N(p-2)}{2}}}{2\|Q\|_2^{p-2}}\|\nabla u\|_2^{\frac{N(p-2)}{2}}+\frac{2}{S^{\frac{2^*}{2}}}\|\nabla u\|_2^{2^*}.
\end{align}

For the case of $p>2+\frac{4}{N}$,
by $\displaystyle\Omega_{a,b,0,\frac{1}{S^{\frac{2^*}{2}}}}^{2,4,p,2^*}>1$ and Lemma \ref{lemma:20220923-hhl1}, there exists $\delta>0$ small such that $\displaystyle\Omega_{a-\delta,b-\delta,0,\frac{1}{S^{\frac{2^*}{2}}}}^{2,4,p,2^*}>1$. So
\begin{align}\label{4.33}
(a-\delta)\|\nabla u\|_2^2 + (b-\delta)\|\nabla u\|_2^4 \ge \frac{1}{S^{\frac{2^*}{2}}}\|\nabla u\|_2^{2^*}.
\end{align}
By (\ref{4.32})-(\ref{4.33}), we obtain that
\begin{align}\lab{eq:20220927-zbe1}
2 \delta\|\nabla u\|_2^2+2 \delta\|\nabla u\|_2^4 \le \frac{N(p-2) c^{p-\frac{N(p-2)}{2}}}{2\|Q\|_2^{p-2}}\|\nabla u\|_2^{\frac{N(p-2)}{2}}.
\end{align}
However, by Lemma \ref{lemma:20220923-hhl1} again, we have that $\displaystyle \Omega_{2\delta, 2\delta,1,0}^{2,4,\frac{N(p-2)}{2},2^*}>0$ and thus
\beq\lab{eq:20220927-zbe2}
2 \delta\|\nabla u\|_2^2+2 \delta\|\nabla u\|_2^4\geq \Omega_{2\delta, 2\delta,1,0}^{2,4,\frac{N(p-2)}{2},2^*} \|\nabla u\|_{2}^{\frac{N(p-2)}{2}}.
\eeq
By $p<2^*$, we have that $p-\frac{N(p-2)}{2}>0$, and thus $\displaystyle \frac{N(p-2) c^{p-\frac{N(p-2)}{2}}}{2\|Q\|_2^{p-2}}< \Omega_{2\delta, 2\delta,1,0}^{2,4,\frac{N(p-2)}{2},2^*}$ for all $\displaystyle c < c_0:=\left(\frac{2 \|Q\|_2^{p-2} \Omega_{2\delta, 2\delta,1,0}^{2,4,\frac{N(p-2)}{2},2^*}}{N(p-2)}\right)^{\frac{2}{2p-N(p-2)}}$. So by \eqref{eq:20220927-zbe1}-\eqref{eq:20220927-zbe2}, we obtain that $\|\nabla u\|_2=0$. Hence, $u\equiv 0$, a contradiction.

For the case of $p=2+\frac{4}{N}$, by \eqref{4.32}, we have that
\beq\lab{eq:20220928-zwe1}
\left(a-\frac{c^{\frac{4}{N}}}{\|Q\|_{2}^{\frac{4}{N}}}\right)\|\nabla u\|_2^2+b\|\nabla u\|_2^4 \leq \frac{1}{S^{\frac{2^*}{2}}}\|\nabla u\|_2^{2^*}.
\eeq
Let $\displaystyle c_0:=\|Q\|_2 \left[a -\frac{4-2^*}{2S^{\frac{N}{N-4}}}\left(\frac{2^*-2}{2b}\right)^{\frac{2}{N-4}}\right]^{\frac{N}{4}}$. By Lemma \ref{lemma:20220923-hhl1}, we get
$$\Omega_{a-\frac{c^{\frac{4}{N}}}{\|Q\|_{2}^{\frac{4}{N}}},b,0,\frac{1}{S^{\frac{2^*}{2}}}}^{2,4,p,2^*}>1, \ \forall c\in (0,c_0).$$
Then by the definition of $\displaystyle \Omega_{a-\frac{c^{\frac{4}{N}}}{\|Q\|_{2}^{\frac{4}{N}}},b,0,\frac{1}{S^{\frac{2^*}{2}}}}^{2,4,p,2^*}$, the formula \eqref{eq:20220928-zwe1} implies that for $c\in (0,c_0)$, $\|\nabla u\|_2=0$ and thus $u\equiv 0$, contradiction.
\ep

\subsection{The case of $p>2+\frac{4}{N}$}
\subsubsection{The normalized ground state solution}
\bl\lab{lemma:20220922-l3} Let $N\geq 5$, $\displaystyle \Omega_{\frac{a}{2},\frac{b}{4},0,\frac{1}{2^* S^{\frac{2^*}{2}}}}^{2,4,\frac{N(p-2)}{2},2^*}>1$ and $p \in \left(2+\frac{4}{N},2^*\right)$. Let $c_1$ be defined by \eqref{def:c1}. Then $c_1>0$.
\el
\bp
We note that $\frac{N(p-2)}{2}>2$ for $p>2+\frac{4}{N}$.
Since $\displaystyle \Omega_{\frac{a}{2},\frac{b}{4},0,\frac{1}{2^* S^{\frac{2^*}{2}}}}^{2,4,\frac{N(p-2)}{2},2^*}>1$, by Lemma \ref{lemma:20220923-hhl1}, there exists $\sigma>0$ such that
 $$\Omega_{\frac{a}{2},\frac{b}{4},\kappa,\frac{1}{2^* S^{\frac{2^*}{2}}}}^{2,4,\frac{N(p-2)}{2},2^*}>1, \forall \kappa\in [0,\sigma).$$
 By $p<2^*$, we see that $p-\frac{N(p-2)}{2}>0$. Let $c_0:=\left(2 \|Q\|_2^{p-2} \sigma\right)^{\frac{2}{2p-N(p-2)}}$. So
 $$\frac{c^{p-\frac{N(p-2)}{2}}}{2 \|Q\|_2^{p-2}}<\sigma, \forall c\in (0,c_0).$$
 Hence,
 $$\Omega_{\frac{a}{2},\frac{b}{4},\kappa,\frac{1}{2^* S^{\frac{2^*}{2}}}}^{2,4,\frac{N(p-2)}{2},2^*}>1, ~\hbox{with}~\kappa=\frac{c^{p-\frac{N(p-2)}{2}}}{2 \|Q\|_2^{p-2}}~\hbox{and}~0<c<c_0.$$
 Then by \eqref{4.2} and the definition of $\Omega_{\frac{a}{2},\frac{b}{4},\kappa,\frac{1}{2^* S^{\frac{2^*}{2}}}}^{2,4,\frac{N(p-2)}{2},2^*}$, we obtain that for $c \in (0,c_0)$, $I(u)>0$ for all $u\in S_c$ and thus $I_c\geq 0$.
 Combining with Lemma \ref{lemma:property-Ic}-(i), we obtain that $I_c=0$ provided $0<c<c_0$. Hence, $c_1\geq c_0>0$.
\ep


Furthermore, we have the following property.
\bt\lab{th:20220923-t1} Let $N\geq 5$, $\displaystyle\Omega_{a,b,0,\frac{1}{S^{\frac{2^*}{2}}}}^{2,4,p,2^*}>1$ and $p \in \left(2+\frac{4}{N},2^*\right)$.
Then
\beq\lab{def:c1-another}
c_1=\sup\{c>0:I_c=0\}.
\eeq
And $I_c=0$ for $c \in (0,c_1]$ while $I_c<0$ for $c>c_1$.
In particular,
\begin{itemize}
\item[(i)]$I_c=0$ and it is not attained provided $0<c<c_1$.
\item[(ii)]$I_c<0$ and it is attained when $c> c_1$.
\item[(iii)] $I_{c_1}=0$ and it is attained.
\end{itemize}
\et

\bp
By \eqref{eq:20220926-e3}, we have that
$\Omega_{\frac{a}{2},\frac{b}{4},0,\frac{1}{2^* S^{\frac{2^*}{2}}}}^{2,4,\frac{N(p-2)}{2},2^*}>\Omega_{a,b,0,\frac{1}{S^{\frac{2^*}{2}}}}^{2,4,p,2^*}>1$.
So by Lemma \ref{lemma:20220927-zbbl1} and Lemma \ref{lemma:20220922-l3}, we see that $c_1\in (0,+\infty)$. By Lemma \ref{lemma:property-Ic}-(ii) and (iii), the monotonicity and continuity imply that
$$I_c\begin{cases}
=0\quad& \hbox{if}~0<c\leq c_1,\\
<0&\hbox{if}~c>c_1.
\end{cases}$$
So \eqref{def:c1-another} holds.
In particular, the conclusion of (i) has been stated out in Lemma \ref{lemma:20220927-zbbl1}.
We also note that (ii) will be proved by Lemma \ref{lemma:20220923-xl1} and (iii) will be proved by Lemma \ref{lemma:20220923-xl2} in the following.
\ep

\bl\lab{lemma:20220923-xl1} Let $N\geq 5$, $\displaystyle\Omega_{a,b,0,\frac{1}{S^{\frac{2^*}{2}}}}^{2,4,p,2^*}>1$ and $p \in \left(2+\frac{4}{N},2^*\right)$. Then $I_c$ is attained by some $u\in S_{r,c}$ if $c>c_1$.
\el

\bp
Noting that we have proved that $-\infty<I_c<0$ for $c>c_1$ in Lemma \ref{lemma:property-Ic} and Theorem \ref{th:20220923-t1}.
Let $\{u_n\} \subset H^1(\R^N)\cap S_c$ be a minimizing sequence for $I_c$. That is, $\|u_n\|_2 =c$ and $I(u_n)=I_c+o_n(1)$. Without loss of generality, we may assume that $u_n \in S_{r,c}$. Indeed,
let $u_n^*$ be the Schwarz rearrangement of $u_n$. Then $\{u_n^*\} \subset H^1_r(\R^N)$ is also a minimizing sequence for $I_c$. Obviously, $I(u_n)<0$ for $n$ large.
Then there exists $t_n \in \R$ such that $T(u_n,t_n) \in P_{r,c}$
and
\begin{align}\label{4.18}
I(T(u_n,t_n))=\min\{I(T(u_n,t)):t \in \R\}<0.
\end{align}
Let $v_n=T(u_n,t_n)$. Then $I(v_n) \le I(u_n)$. Moreover, $I(v_n) \rightarrow I_c$.
By the Ekeland principle, there exists $\{w_n\} \subset S_{r,c}$ such that
$I(w_n) \le I(v_n)$, $\|w_n-v_n\|_{H^1} \rightarrow 0$ and for all $v \in S_{r,c}$,
\begin{align}\label{eq:20220923-hwbe2}
I(w_n) \le I(v)+ \frac{1}{n} \|w_n-v\|_{H^1}.
\end{align}
It follows from the implicit function theorem that $S_{r,c}$ is a $C^1$-manifold of codimension 1 and $H_r^1(\R^N)=\R u\oplus T_{u}(S_{r,c})$ for each $u\in S_{r,c}$. Fix $u\in S_{r,c}$, then for any $\phi\in H_r^1(\R^N)$, we write $\phi=\psi+s_{u,\phi} u$ with $\psi\in T_u(S_{r,c})$ and $s_{u,\phi}\in \R$. Precisely, $s_{u,\phi}$ can be computed by
\beq\lab{eq:20220923-wbe1}
s_{u,\phi}=\frac{\int_{\R^N}\phi u dx}{\|u\|_2^2}=\frac{1}{c^2}\int_{\R^N}\phi u dx.
\eeq
Now, we fix $u=w_n$ and consider $v(t)=c \frac{w_n+ t \phi}{\|w_n+ t \phi\|_2}$, $t\in \R$. Then $v(t)$ is a $C^1$ curve in $S_{r,c}$. In particular, $v(0) =w_n$ and
$v'(0)=\phi-\frac{\int_{\R^N}w_n \phi \mathrm{d}x }{c^2}w_n=\phi-s_{w_n,\phi}w_n\in T_{w_n}(S_{r,c})$. By a direct calculation, the formula \eqref{eq:20220923-hwbe2} implies that
\begin{align}\label{4.20}
&\left(I'(w_n),\phi\right)-\frac{1}{c^2}\left(I'(w_n),w_n\right)\int_{\R^N}w_n \phi \mathrm{d}x\notag\\
& \ge -\frac{1}{n} \left\|\phi-\frac{\int_{\R^N}w_n \phi \mathrm{d}x}{c^2}w_n\right\|_{H^1},\forall \phi\in H_r^1(\R^N).
\end{align}
Thus
\begin{align}\label{eq:20220923-we1}
&\left(I'(w_n),\phi-s_{w_n,\phi}w_n\right)-\frac{1}{c^2}\left(I'(w_n),w_n\right)\int_{\R^N}w_n (\phi-s_{w_n,\phi}w_n) \mathrm{d}x\notag\\
& \ge -\frac{1}{n} \left\|\phi-s_{w_n,\phi}w_n\right\|_{H^1}.
\end{align}
By the arbitrary of $\phi\in H_r^1(\R^N)$, we obtain the arbitrary of $\psi=\phi-s_{w_n,\phi}w_n\in T_{w_n}(S_{r,c})$.
Hence, there exists $\{\lambda_n\} \subset \R$ with $\lambda_n:=\frac{1}{c^2}(I'(w_n), w_n)$ such that
\begin{align}\label{eq:20220923-wbe2}
\|I'(w_n) - \lambda_n \varphi'(w_n)\|_{(T_{w_n}S_{r,c})^*} \rightarrow 0,
\end{align}
where $\varphi(w_n)=\frac{1}{2}\|w_n\|_2^2$.
By the definition of $\lambda_n$, we have that $(I'(w_n) - \lambda_n \varphi'(w_n), w_n)\equiv 0$. So by $H_r^1(\R^N)=\R w_n\oplus T_{w_n}(S_{r,c})$, the formula \eqref{eq:20220923-wbe2} implies that
$\displaystyle I'(w_n) - \lambda_n \varphi'(w_n)\rightarrow 0$ as $n\rightarrow +\infty$.
That is, $I|'_{S_{r,c}}(w_n)\rightarrow 0$ as $n\rightarrow +\infty$.

Since $v_n \in P_{r,c}$, by $\|w_n-v_n\|_{H^1} \rightarrow 0$, we have that $G(w_n) \rightarrow 0$.
Together with $I(w_n) \rightarrow I_c$ and $I|'_{S_{r,c}}(w_n)\rightarrow 0$, we see that $\{w_n\}\subset S_{r,c}$ is a $(PSP)_{I_c}$ sequence. Since $I_c<0$ for $c>c_1$, by Lemma \ref{lemma:20220923-zl1}, $\{w_n\}$ converges strongly to some $u$ in $H_r^1(\R^N)$ up to a subsequence. Hence,
$I_c$ is attained by $u\in S_{r,c}$.
\ep

\bl\lab{lemma:20220923-xl2} Let $N\geq 5$, $\displaystyle\Omega_{a,b,0,\frac{1}{S^{\frac{2^*}{2}}}}^{2,4,p,2^*}>1$ and $p \in \left(2+\frac{4}{N},2^*\right)$. Then $I_{c_1}=0$ is attained by some $u\in S_{r,c_1}$.
\el

\bp Let $c_n=c_1+\frac{1}{n}$. By Lemma \ref{lemma:20220923-xl1}, there exists $\{u_n\} \subset S_{r,c_n}$ such that $I(u_n)=I_{c_n}<0$.  By Lemma \ref{lemma:property-Ic} and Theorem \ref{th:20220923-t1}, we have that $\displaystyle\lim_{n \rightarrow \infty}I(u_n)=\lim_{n \rightarrow \infty}I_{c_n}=I_{c_1}=0$.
By \eqref{eq:20220923-zbe1}, we obtain that $\{u_n\}$ is bounded in $H_r^1(\R^N)$. Up to a subsequence, we assume that $u_n \rightharpoonup u_0$ weakly in $H_r^1(\R^N)$. If $u_0=0$, we have that $\displaystyle\lim_{n \rightarrow \infty}\|u_n\|_p = 0$. Thanks to $\displaystyle\Omega_{a,b,0,\frac{1}{S^{\frac{2^*}{2}}}}^{2,4,p,2^*}>1$, by a similar argument as in the proof of Lemma \ref{lemma:20220923-zl1}, we can get that
$\displaystyle\lim_{n \rightarrow \infty}\|\nabla u_n\|_2 = 0$.
However, by $p>2+\frac{4}{N}$, we have that $\frac{N(p-2)}{2}>2$. Then by \eqref{4.2}, we have that $I_{c_n}=I(u_n)>0$ for $n$ large enough, a contradiction.
So $u_0 \ne 0$. Since $I_{c_n}$ is attained by $u_n$, there exists $\lambda_n \in \R$ such that in a weak sense,
\begin{align}\label{4.22}
-\left(a+b\|\nabla u_n\|_2^2\right)\Delta u_n=\lambda_n u_n+ |u_n|^{p-2}u_n+|u_n|^{2^*-2}u_n~\hbox{in}~\R^N.
\end{align}
Applying a similar argument in the proof of Lemma \ref{lemma:20220923-zl1}, we can have that
\beq
\lambda_n c_n^2 =\left[\frac{N}{2}-\frac{N}{p}-1\right]\|u_n\|_p^p\rightarrow \left[\frac{N}{2}-\frac{N}{p}-1\right]\|u_0\|_p^p.
\eeq
Since $\lim_{n\rightarrow +\infty}c_n=c_1>0$, we get
\beq
\lambda_n\rightarrow \frac{1}{c_1^2}\left[\frac{N}{2}-\frac{N}{p}-1\right]\|u_0\|_p^p:=\lambda_0<0.
\eeq
Let $E=\lim_{n \rightarrow \infty}\|
\nabla u_n\|_2^2$. Then $u_0$ weakly solves \eqref{eq:20220923-zbe2} and the equalities \eqref{eq:20220923-zbe3} and \eqref{eq:20220923-zbe4} hold.
Let $v_n:=u_n-u_0$, we can deduce \eqref{4.17} again.
Then $\|\nabla v_n\|_2\rightarrow 0$, which means that $u_n\rightarrow u_0$ in $D^{1,2}(\R^N)$. Hence, $E=\|\nabla u_0\|_2^2$. Then one can see that
$$I(u_0)=\lim_{n\rightarrow +\infty}I(u_n)=I_{c_1}.$$
In particular, similar to \eqref{eq:20220923-zbe8}, we can prove that
$$\lambda_0(\|u_0\|_2^2-c_1^2)=0.$$
And thus $\|u_0\|_2^2=c_1^2$. That is, $u_0\in S_{r,c_1}$ attains $I_{c_1}=0$.
\ep

\subsubsection{The mountain pass type solution}

\bl\lab{lemma:20220923-wbbl1} Let $N\geq 5$, $\displaystyle\Omega_{a,b,0,\frac{1}{S^{\frac{2^*}{2}}}}^{2,4,p,2^*}>1$ and $p \in \left(2+\frac{4}{N},2^*\right)$. Then there exist $\theta>1$ and $\eta \in (0,1)$ small such that for any $K_c \in (0,\eta)$, it holds that
\begin{align}\label{4.26}
0< \sup_{u \in A_c}I(u) < \inf_{u \in B_c}I(u),
\end{align}
where
\begin{align*}
A_c:=\{u \in S_{r,c}: \|\nabla u\|_2^2 \le K_c\},\ \
B_c:=\{u \in S_{r,c}: \|\nabla u\|_2^2 = \theta K_c\}.
\end{align*}
Moreover, $I(u)>0$ for all $u \in A_c$.
\el

\bp By \eqref{eq:20220926-e3}, we have that
$\displaystyle\Omega_{\frac{a}{2},\frac{b}{4},0,\frac{1}{2^*S^{\frac{2^*}{2}}}}^{2,4,p,2^*}>\Omega_{a,b,0,\frac{1}{S^{\frac{2^*}{2}}}}^{2,4,p,2^*}>1$.
By Lemma \ref{lemma:20220923-hhl1}, there exists some $\delta\in (0,1)$ small such that
$\displaystyle \Omega_{\frac{a}{2}(1-\delta),\frac{b}{4}(1-\delta^2),0,\frac{1}{2^*S^{\frac{2^*}{2}}}}^{2,4,p,2^*}>1$.
So
\beq\lab{eq:20220926-xe1}
\frac{a}{2}(1-\delta)t^2+\frac{b}{4}(1-\delta^2) t^4>\frac{1}{2^*S^{\frac{2^*}{2}}} t^{2^*}, \ \forall t>0.
\eeq
For any $u\in A_c$, we have that
\begin{align}\label{347}
I(u)
\geq &\frac{a}{2}\|\nabla u\|_2^2+\frac{b}{4}\|\nabla u\|_2^4-\frac{1}{p}\|u\|_p^p -\frac{1}{2^*S^{\frac{2^*}{2}}}\|\nabla u\|_{2}^{2^*}\quad \hbox{by \eqref{4.1}}\notag\\
>&\frac{a \delta}{2} \|\nabla u\|_2^2+\frac{b \delta^2}{4}\|\nabla u\|_2^4-\frac{1}{p}\|u\|_p^p\quad \hbox{by \eqref{eq:20220926-xe1}}\notag\\
\geq&\frac{a \delta}{2} \|\nabla u\|_2^2+\frac{b \delta^2}{4}\|\nabla u\|_2^4-\frac{1}{2\|Q\|_{2}^{p-2}}\|\nabla u\|_{2}^{\frac{N(p-2)}{2}} c^{p-\frac{N(p-2)}{2}}\quad \hbox{by Lemma \ref{lemma:20220923-l1}}.
\end{align}
Since $p > 2+\frac{4}{N}$, we see that $\frac{N(p-2)}{2}>2$. Hence, there exists $\eta_1>0$ small such that for $K_c\in (0,\eta_1)$, we have that $I(u)>0$ for any $u\in A_c$.

Furthermore, we take $\theta=\frac{2}{\delta}$, then $\theta>1$ and for any $v\in B_c$ and $ u\in A_c$,
\begin{align*}
I(v)-I(u)=&\frac{a}{2} \left(\|\nabla v\|_2^2-\|\nabla u\|_2^2\right)+ \frac{b}{4} \left(\|\nabla v\|_2^4-\|\nabla u\|_2^4\right)\\
&-\frac{1}{2^*}\left(\|v\|_{2^*}^{2^*}-\|u\|_{2^*}^{2^*}\right)-\frac{1}{p}\left(\|v\|_p^p-\|u\|_p^p\right)\\
>&\frac{a}{2} \left(\|\nabla v\|_2^2-\|\nabla u\|_2^2\right)+ \frac{b}{4} \left(\|\nabla v\|_2^4-\|\nabla u\|_2^4\right)-\frac{1}{2^*}\|v\|_{2^*}^{2^*}-\frac{1}{p}\|v\|_p^p\\
>&\frac{a \delta}{2} \|\nabla v\|_2^2+\frac{b \delta^2}{4} \|\nabla v\|_2^4 -\frac{a}{2}\|\nabla u\|_2^2-\frac{b}{4}\|\nabla u\|_2^4-\frac{1}{p}\|v\|_p^p\\
\geq&\frac{a \delta}{2}\|\nabla v\|_2^2+\frac{b \delta^2}{4} \|\nabla v\|_2^4
-\frac{a K_c}{2}-\frac{b K_c^2}{4}-\frac{1}{2\|Q\|_{2}^{p-2}}\|\nabla v\|_{2}^{\frac{N(p-2)}{2}} c^{p-\frac{N(p-2)}{2}}\\
=&\frac{aK_c}{2}+\frac{3bK_c^2}{4}-\frac{1}{2\|Q\|_{2}^{p-2}} c^{p-\frac{N(p-2)}{2}} \left(\frac{2}{\delta}\right)^{\frac{N(p-2)}{4}} K_{c}^{\frac{N(p-2)}{4}}\\
=:&\Lambda(K_c,\delta)
\end{align*}
Then by $\frac{N(p-2)}{4}>1$, there exists $\eta_2\in (0,1)$ small such that for $K_c\in (0,\eta_2)$, we have that
$$I(v)-I(u)>\Lambda(K_c,\delta)>0,\  \forall v\in B_c,\ \forall u\in A_c.$$
And thus
$\displaystyle\sup_{u\in A_c}I(u)<\inf_{u\in B_c}I(u)$.
Hence, we can take $\eta:=\min\{\eta_1,\eta_2\}$ and the proof is finished.
\ep

\br\lab{remark:20220927-r1}
By the details in the proof of Lemma \ref{lemma:20220923-wbbl1}, the number $\eta$ indeed can be chosen uniformly for bounded $c$. That is, for such a fixed $\theta>1$, we can find $\eta:=\eta_M>0$ small such that
\beq\lab{eq:20220927-e1}
0<\sup_{u \in S_{r,c}, \|\nabla u\|_2^2 \le K_c}I(u)<\inf_{u \in S_{r,c}, \|\nabla u\|_2^2 = \theta K_c}I(u), \ \forall c\leq M, \ \forall K_c \leq \eta.
\eeq
\er

\bt\lab{th:20220926-bt1} Let $N\geq 5$, $\displaystyle\Omega_{a,b,0,\frac{1}{S^{\frac{2^*}{2}}}}^{2,4,p,2^*}>1$ and $p \in \left(2+\frac{4}{N},2^*\right)$. Then $I|_{S_{r,c}}$ has a mountain pass type critical point for any $c\geq c_1$.
\et

\bp
By Theorem \ref{th:20220923-t1}, for any $c\geq c_1$, there exists $u_0 \in S_{r,c}$ such that $I(u_0)=I_c \leq 0$. Let $\theta>1$ and $\eta>0$ be given by Lemma \ref{lemma:20220923-wbbl1},
we take
$$K_c<\min\left\{\frac{1}{\theta}\|\nabla u_0\|_2^2, \eta\right\}.$$
Then $\|\nabla u_0\|_2^2 > \theta K_c$. Noting that $\displaystyle\lim_{s \rightarrow -\infty}\|\nabla T(u_0,s)\|_2=0$. We can take some $s_1<0$ and put $u_1:=T(u_0,s_1)$ such that $u_1 \in S_{r,c}$ and $\|\nabla u_1\|_2^2 \le \frac{1}{2}K_c$.
Define
\begin{align}\label{4.28}
\gamma_c:=\inf_{\gamma \in \Gamma_c}\max_{t \in [0,1]}I(\gamma(t)),
\end{align}
where $\Gamma_c:=\left\{\gamma \in C([0,1],S_{r,c}):\gamma(0)=u_1, \gamma(1)=u_0\right\}$.
By the choice of $K_c$ and Lemma \ref{lemma:20220923-wbbl1}, we have that
\beq
\gamma_c>\max\{I(u_1), I(u_0)\}.
\eeq
Similar to the argument of \cite[Proposition 2.2 and Lemma 2.4]{J}, we can derive the existence of $(PSP)_{\gamma_c}$-sequence, that is, $\{u_n\}\subset S_{r,c}$ satisfies
 $$I(u_n)
\rightarrow \gamma_c, I|'_{S_{r,c}}(u_n) \rightarrow 0~\hbox{and}~G(u_n) \rightarrow 0.$$
Noting that $\gamma_c>0$, by Lemma \ref{lemma:20220923-zl1}, there exists $u \in H_r^1(\R^N)$ such that, up to a subsequence, $u_n \rightarrow u$ in $H_r^1(\R^N)$. So $u\in S_{r,c}$, $I(u)
= \gamma_c$ and $I|'_{S_{r,c}}(u)=0$.
\ep

\bt\lab{th:20220926-bt2} Let $N\geq 5$, $\displaystyle\Omega_{a,b,0,\frac{1}{S^{\frac{2^*}{2}}}}^{2,4,p,2^*}>1$ and $p \in \left(2+\frac{4}{N},2^*\right)$. Then there exists $\bar{\eta}\in (0,1)$ small such that $I|_{S_{r,c}}$ has a mountain pass type critical point for $c \in [c_1-\bar{\eta},c_1)$.
\et
\bp
By Remark \ref{remark:20220927-r1}, there exist $\theta>1$ and $\eta \in (0,1)$ small independent of $c\leq c_1$ such that
\beq\lab{eq:20220927-e2}
0<\sup_{u\in S_{r,c}, \|\nabla u\|_2^2\leq K_c}I(u)<\inf_{u\in S_{r,c}, \|\nabla u\|_2^2=\theta K_c}I(u), \ \forall c\leq c_1, \ \forall K_c\leq \eta.
\eeq
By Lemma \ref{lemma:20220923-xl2}, there exists some $u_0\in S_{r,c_1}$ such that $I(u_0)=I_{c_1}=0$.
By (\ref{347}), we can choose $K_{c_1} \in (0,\eta]$ small and $\varepsilon_0 \in (0,1)$ small such that $\|\nabla u_0\|_2^2 >\theta  K_{c_1}$ and
\beq\lab{eq:20220927-zwe2}
\inf_{u\in S_{r, (1-\varepsilon)c_1}, \|\nabla u\|_2^2=\theta (1-\varepsilon)^2 K_{c_1}}I(u)>I((1-\varepsilon)u_0), \ \forall \varepsilon\in [0,\varepsilon_0].
\eeq
For $c<c_1$ close to $c_1$, we rewrite $c=(1-\varepsilon)c_1$ with $\varepsilon>0$ small.
By \eqref{eq:20220927-e2}, we can take $K_c:=(1-\varepsilon)^2 K_{c_1}<K_{c_1}\leq \eta$ such that
\beq\lab{eq:20220927-e3}
0<\sup_{u\in S_{r,c},\|\nabla u\|_2^2\leq  K_c}I(u)<\inf_{u\in S_{r,c},\|\nabla u\|_2^2=\theta K_c}I(u).
\eeq
Put $u_1:=(1-\varepsilon)u_0\in S_{r,c}$, we note that $\|\nabla u_1\|_2^2=(1-\varepsilon)^2 \|\nabla u_0\|_2^2>(1-\varepsilon)^2 \theta K_{c_1}=\theta K_c$.
Let $u_2:=T(u_1, s_2)$ with some $s_2<0$  such that  $\|\nabla u_2\|_2^2 \le \frac{1}{2}K_c$. Then for $\varepsilon\in [0,\varepsilon_0]$, it holds that
\beq\lab{eq:20220927-zwe1}
\inf_{u\in S_{r,c},\|\nabla u\|_2^2=\theta K_c}I(u)>\max\{I(u_1), I(u_2)\}.
\eeq
Hence, we have that
\beq
\gamma_c:=\inf_{\gamma\in \Gamma_c}\max_{t\in [0,1]}I(\gamma(t))\geq \inf_{u\in S_{r,c},\|\nabla u\|_2^2=\theta K_c}I(u)>\max\{I(u_1), I(u_2)\},
\eeq
where $\displaystyle\Gamma_c:=\left\{\gamma\in C([0,1], S_{r,c}): \gamma(0)=u_1, \gamma(1)=u_0\right\}$.
Similar to the argument of Theorem \ref{th:20220926-bt1}, we obtain the final result.
\ep

\subsubsection{The local constraint minimizer}

\bt\lab{lemma:221225} Let $N\geq 5$, $\displaystyle\Omega_{a,b,0,\frac{1}{S^{\frac{2^*}{2}}}}^{2,4,p,2^*}>1$ and $p \in \left(2+\frac{4}{N},2^*\right)$. Then
$I|_{S_{r,c}}$ has a local constraint minimizer with positive energy for $c \in [c_1-\bar{\eta},c_1)$.
\et

\bp By Lemma \ref{lemma:20220923-xl2}, there exists some $u_0\in S_{r,c_1}$ such that $I(u_0)=I_{c_1}=0$.
By (\ref{eq:20220926-e3}), (\ref{347}) and $\displaystyle\Omega_{a,b,0,\frac{1}{S^{\frac{2^*}{2}}}}^{2,4,p,2^*}>1$, we can choose $K_{c_1} \in (0,\eta]$ small and $\varepsilon_0 \in (0,\frac{1}{2})$ small such that $\|\nabla u_0\|_2^2 >\theta  K_{c_1}$,
\beq\lab{eq:1226-1}
\inf_{u\in S_{r, (1-\varepsilon)c_1}, \|\nabla u\|_2^2=\theta (1-\varepsilon)^2 K_{c_1}}I(u)>I((1-\varepsilon)u_0), \ \forall \varepsilon\in [0,\varepsilon_0],
\eeq
and
\beq\lab{eq:1226-2}
\inf_{\|\nabla u\|_2^2\ge \theta (1-\varepsilon)^2 K_{c_1}}\left(\frac{a}{2}\|\nabla u\|_2^2+\frac{b}{4}\|\nabla u\|_2^4-\frac{1}{2^* S^{\frac{2^*}{2}}}\|\nabla u\|_2^{2^*}\right)>I((1-\varepsilon)u_0), \ \forall \varepsilon\in [0,\varepsilon_0].
\eeq
Let $c:=(1-\varepsilon)c_1$, $u_1:=(1-\varepsilon)u_0\in S_{r,c}$ and $K_c:=(1-\varepsilon)^2 K_{c_1}$. Define
\begin{align*}
m_c:=\inf_{u\in S_{r,c},\|\nabla u\|_2^2 \ge \theta K_c}I(u).
\end{align*}
By Theorem \ref{th:20220923-t1} and \eqref{eq:1226-1}, we get
\begin{align*}
0 \le m_c \le I(u_1) <\inf_{u\in S_{r,c},\|\nabla u\|_2^2 = \theta K_c}I(u).
\end{align*}
Moreover, there exists $\sigma \in (0,1)$ small such that
\begin{align*}
I(u_1) <\inf_{u\in S_{r,c}, \theta K_c \le \|\nabla u\|_2^2 \le \theta K_c+  \sigma}I(u).
\end{align*}
By the definition of $m_c$, there exists $\{u_n\} \subset S_{r,c}$ such that $\|\nabla u_n\|_2^2 >  \theta K_c+ \sigma$ and $I(u_n) \rightarrow m_c$.
Similar to the argument of Lemma \ref{lemma:20220923-xl1}, there exists $\{w_n\} \subset S_{r,c}$ such that
$\|w_n-u_n\|_{H^1} \rightarrow 0$, $I(w_n) \rightarrow m_c$ and $I|'_{S_{r,c}}(w_n)\rightarrow 0$.

If $w_n \rightharpoonup 0$ weakly in $H_r^1(\R^N)$, then
\begin{align*}
I(u_1) \ge m_c \ge &\frac{a}{2}\lim_{n \rightarrow \infty}\|\nabla w_n\|_2^2+\frac{b}{4}\lim_{n \rightarrow \infty}\|\nabla w_n\|_2^4-\frac{1}{2^* S^{\frac{2^*}{2}}}\lim_{n \rightarrow \infty}\|\nabla w_n\|_2^{2^*},
\end{align*}
a contradiction with $\lim_{n \rightarrow \infty}\|\nabla w_n\|_2^2 \ge  \theta K_c$ and \eqref{eq:1226-2}. So $w_n \rightharpoonup w \ne 0$ weakly in $H_r^1(\R^N)$. We assume that $\displaystyle E=\lim_{n \rightarrow \infty}\|
\nabla w_n\|_2^2$. By $I|'_{S_{r,c}}(w_n) \rightarrow 0$, we get
 \begin{align*}
 I'(w_n)w_n=\lambda_n\|w_n\|_2^2+o_n(1).
 \end{align*}
 So $|\lambda_n|$ is bounded. Furthermore,
\beq\label{1226-zbe1}
\begin{cases}
&I'(w_n) - \lambda_n w_n \rightarrow 0,\\
&\lambda_n \rightarrow \lambda_0,\\
&\|w_n\|_2^2=c, \|\nabla w_n\|_2^2\rightarrow E,\\
&w_n\rightharpoonup w~\hbox{weakly in $H_r^1(\R^N)$ and $w_n\rightarrow w$ in $L^p(\R^N)$}.
\end{cases}
\eeq
Then $w \in H_r^1(\R^N)$ weakly solves
\beq\lab{1226-zbe2}
-(a+bE)\Delta w-\lambda_0 w=|w|^{p-2}w+|w|^{2^*-2}w~\hbox{in}~\R^N,
\eeq
from which we derive that
\beq\lab{1226-zbe3}
(a+bE)\|\nabla w\|_2^2-\lambda_0\|w\|_2^2=\|w\|_p^p+\|w\|_{2^*}^{2^*},
\eeq
and
\beq\lab{1226-zbe4}
(N-2)(a+bE)\|\nabla w\|_2^2-N\lambda_0\|w\|_2^2=2N\left[\frac{1}{p}\|w\|_p^p+\frac{1}{2^*}\|w\|_{2^*}^{2^*}\right].
\eeq
By \eqref{1226-zbe3}-\eqref{1226-zbe4}, we obtain that
\beq\lab{1226-zbe5}
\lambda_0\|w\|_2^2=\left[\frac{N}{2}-\frac{N}{p}-1\right]\|w\|_p^p.
\eeq
So $\lambda_0<0$. Let $v_n=w_n-w$. By \eqref{1226-zbe1}-\eqref{1226-zbe2}, we have that
\beq\lab{1226-zbe6}
(a+b\|\nabla v_n\|_2^2)\|\nabla v_n\|_2^2-\lambda_0\|v_n\|_2^2\le \|v_n\|_{2^*}^{2^*}+o_n(1) \le \frac{\|\nabla v_n\|_{2}^{2^*}}{S^{\frac{2^*}{2}}}+ o_n(1).
\eeq
So by $\displaystyle\Omega_{a,b,0,S^{-\frac{2^*}{2}}}^{2,4,p,2^*}>1$, we obtain that
 $\|\nabla v_n\|_2\rightarrow 0$. Moreover, $\|v_n\|_2\rightarrow 0$.
Hence, we have that $w_n \rightarrow w$ in $H_r^1(\R^N)$.

By $I(w_n) \rightarrow m_c$, $I|'_{S_{r,c}}(w_n)\rightarrow 0$ and $w_n \rightarrow w$ in $H_r^1(\R^N)$, we get $I(w) = m_c$ and $I|'_{S_{r,c}}(w)= 0$.
By Theorem \ref{th:20220923-t1}, we know that $I_c=0$ is not attained for $c\in (0,c_1)$. Then $m_c>0$ for $c \in [c_1-\bar{\eta},c_1)$.
\ep

\noindent\textit{\bf Proof of  Theorem  \ref{th:main-zzzt1}-(i).} By Theorem \ref{th:20220923-t1}, we get (i-1); by Theorem \ref{th:20220926-bt1}, we get (i-2); by Theorem \ref{th:20220926-bt2} and Theorem \ref{lemma:221225}, we get (i-3); by Theorem \ref{th:20220928-t1}, we get (i-4).
\hfill $\square$

\subsection{The case of $p=2+\frac{4}{N}$}
Put
\begin{align}\label{eq:20220928-ze1}
c_*:=\|Q\|_2 \left[a -\frac{(4-2^*)}{(2^*)^{\frac{N-2}{N-4}}S^{\frac{N}{N-4}}}\left(\frac{2(2^*-2)}{b}\right)^{\frac{2}{N-4}}\right]^{\frac{N}{4}}.
\end{align}
We note that $\frac{2^{\frac{2}{N-4}}}{(2^*)^{\frac{N-2}{N-4}}}<\left(\frac{1}{2}\right)^{\frac{N-2}{N-4}}$ for $N\geq 5$. Suppose that $\displaystyle\Omega_{a,b,0,\frac{1}{S^{\frac{2^*}{2}}}}^{2,4,p,2^*}>1$, then by \eqref{eq:20220924-ze1} in Remark \ref{remark:20220923-wr1}, we can have that
\begin{align*}
&\Omega_{a,b,0,\frac{1}{S^{\frac{2^*}{2}}}}^{2,4,p,2^*}>1\\
\Leftrightarrow&\left(\frac{2a}{4-2^*}\right)^{\frac{4-2^*}{2}} \left(\frac{2b}{2^*-2}\right)^{\frac{2^*-2}{2}}>\frac{1}{S^{\frac{2^*}{2}}}\\
\Leftrightarrow&a>\left(\frac{1}{2}\right)^{\frac{N-2}{N-4}}(4-2^*)\left(\frac{2^*-2}{b}\right)^{\frac{2}{N-4}}\frac{1}{S^{\frac{N}{N-4}}}\\
\Rightarrow& a>\frac{(4-2^*)}{(2^*)^{\frac{N-2}{N-4}}S^{\frac{N}{N-4}}}\left(\frac{2(2^*-2)}{b}\right)^{\frac{2}{N-4}}\\
\Rightarrow&c_*>0.
\end{align*}

\bl\lab{lemma:20220928-l1}
Let $N\geq 5$, $\displaystyle\Omega_{a,b,0,\frac{1}{S^{\frac{2^*}{2}}}}^{2,4,p,2^*}>1$ and $p =2+\frac{4}{N}$. Let $c_1$ be defined by \eqref{def:c1}. It holds that
 \beq\label{eq:20220928-ze2}
 0<c_*\leq c_1\leq a^{\frac{4}{N}}\|Q\|_2.
 \eeq
\el
\bp
Recalling \eqref{4.2}, by $p=2+\frac{4}{N}$, we obtain that
\beq\label{eq:20220928-ze3}
I(u)\geq \left(\frac{a}{2}-\frac{c^{\frac{4}{N}}}{2 \|Q\|_2^{\frac{4}{N}}}\right)\|\nabla u\|_2^2 + \frac{b}{4} \|\nabla u\|_2^4-\frac{1}{2^* S^{\frac{2^*}{2}}}\|\nabla u\|_2^{2^*}
\eeq
Let
$$A:=\frac{a}{2}-\frac{c^{\frac{4}{N}}}{2 \|Q\|_2^{\frac{4}{N}}}>0, \ B:=\frac{b}{4}.$$
Consider the equation
$$\Omega_{A,B,0,\frac{1}{2^* S^{\frac{2^*}{2}}}}^{2,4,p,2^*}=2^* S^{\frac{2^*}{2}} \Omega_{A,B,0,1}^{2,4,p,2^*}=1.$$
By \eqref{eq:20220926-e1}, we have that
\beq\label{eq:20220928-ze4}
2^*S^{\frac{2^*}{2}} 2 (2^*-2)^{-\frac{2^*-2}{2}} (4-2^*)^{-\frac{4-2^*}{2}} \left(\frac{a}{2}-\frac{c^{\frac{4}{N}}}{2 \|Q\|_2^{\frac{4}{N}}}\right)^{\frac{4-2^*}{2}} \left(\frac{b}{4}\right)^{\frac{2^*-2}{2}}=1.
\eeq
By solving \eqref{eq:20220928-ze4}, we obtain that $c=c_*$. Hence, for $c\leq c_*$, we have that $\displaystyle\Omega_{A,B,0,\frac{1}{2^* S^{\frac{2^*}{2}}}}^{2,4,p,2^*}\geq 1$. So  by \eqref{eq:20220928-ze3} and Lemma \ref{lemma:property-Ic}-(i), we obtain that $I_c=0, \forall c\in (0,c_*]$. Then by the definition of $c_1$, we have that $c_1\geq c_*>0$.

On the other hand, by \eqref{4.9},
\begin{align*}
I(Q_t) \leq \frac{a}{2} c^2 t^2 + \frac{b}{4} c^4 t^4-\frac{c^p t^{\frac{Np}{2}-N}}{2 \|Q\|_2^{p-2}}=\left[a-\frac{c^{\frac{4}{N}} }{ \|Q\|_2^{\frac{4}{N}}}\right]\frac{1}{2}c^2t^2+ \frac{b}{4} c^4 t^4.
\end{align*}
So if $c>a^{\frac{4}{N}}\|Q\|_2$, it holds that $a-\frac{c^{\frac{4}{N}} }{\|Q\|_2^{\frac{4}{N}}}<0$. Hence, $\min_{t \ge 0}I(Q_t)<0$, which implies that
$I_c<0$ for $c>a^{\frac{N}{4}} \|Q\|_2$. By the definition of $c_1$, we have that $c_1\leq a^{\frac{4}{N}}\|Q\|_2$.
\ep

\noindent\textit{\bf Proof of  Theorem  \ref{th:main-zzzt1}-(ii).} By Lemma \ref{lemma:20220928-l1}, we see that $c_1>0$ and thus $I_c<0$ for $c>c_1$. By Lemma \ref{lemma:20220927-zbbl1}, we have that $I_c=0$ for $c \in (0,c_1]$ and $I_c$ is not attained for $c \in (0,c_1)$.
Applying a similar argument as the proof of Lemma \ref{lemma:20220923-xl1}, we can get (ii-1); by Theorem \ref{th:20220928-t1}, we get (ii-2).
\hfill $\square$

\subsection{The case of $2<p<2+\frac{4}{N}$}
By Lemma \ref{lemma:property-Ic}-(i), we get $I_c>-\infty$ for any $c>0$. Furthermore, for the case of $p<2+\frac{4}{N}$,  by \eqref{4.9},
\begin{align*}
I(Q_t) \le \frac{a}{2} c^2 t^2 + \frac{b}{4} c^4 t^4-\frac{c^p t^{\frac{Np}{2}-N}}{2 \|Q\|_2^{p-2}}.
\end{align*}
Then
$$I_c\leq \min_{t>0}I(Q_t)<0.$$

\noindent\textit{\bf Proof of  Theorem  \ref{th:main-zzzt1}-(iii).} In such a case, it holds that $I_c \in (-\infty,0)$ for any $c>0$. Then apply a similar argument as that in Lemma \ref{lemma:20220923-xl1}, we can prove that $I_c$ is attained.
\hfill $\square$

\s{The case $N=4$}\lab{sec:caseN=4}

\renewcommand{\theequation}{4.\arabic{equation}}

\subsection{Some preliminaries}
When $N=4$, by (\ref{4.2}),
\begin{align}\label{eq:20220928-xe1}
I(u) \ge \frac{a}{2}\|\nabla u\|_2^2 + \frac{1}{4} \left(b-\frac{1}{S^2}\right)\|\nabla u\|_2^4-\frac{c^{4-p}}{2 \|Q\|_2^{p-2}}\|\nabla u\|_2^{2(p-2)}
\end{align}
holds for all $u \in S_c$.

If $b>\frac{1}{S^2}$ and $2<p<4=2^*$, similar to the arguments in subsection \ref{subsection:N-5}, we have the following results. Since the proofs are almost the same, we omit the details.

\bl\lab{lemma:20220928-zwl1}
Let $N=4$, $b>\frac{1}{S^2}$ and $p \in (2,4)$. Then
\begin{itemize}
\item[(i)]$-\infty<I_c\leq 0$ for all $c>0$;
\item[(ii)] $I_c$ is non-increasing on $c \in (0,+\infty)$;
\item[(iii)]$I_c$ is continuous on $c \in (0,+\infty)$.
\end{itemize}
\el

\bc\lab{cor:20220928-c1}
Let $N=4$, $b>\frac{1}{S^2}$ and $p \in (2,4)$.
If there exists $0<\underbar{c}<\bar{c}<+\infty$ such that $I_{\underbar{c}}=I_{\bar{c}}$, then $I_c$ is not attained provided $c\in [\underbar{c}, \bar{c})$.
\ec

\bl\lab{lemma:20220928-zwl2}
Let $N=4$, $b>\frac{1}{S^2}$ and $p \in (2,4)$. Then $I_c<0$ for $c>0$ large enough. And thus we can define $c_1$ as in \eqref{def:c1}.
Furthermore, if $c_1>0$, then $I_c\equiv 0$ for $c\in (0,c_1]$ and $I_c$ is not attained provided $c\in (0,c_1)$.
\el

Similar to Lemma \ref{lemma:20220923-zl1}, we also can establish the following compactness result.
\bl\lab{lemma:20220928-zwl3} Let $N=4$, $b>\frac{1}{S^2}$ and $p \in (2,4)$. If $\{u_{n}\}
\subset S_{r,c}$ such that $I(u_n)
\rightarrow \gamma_c \ne 0$, $I'_{S_{r,c}}(u_n) \rightarrow 0$ and $G(u_n) \rightarrow 0$, then $\{u_{n}\}$ converges strongly in $H_r^1(\R^4)$ up to a subsequence.
\el

\subsection{The case of $3<p<4$}

\bl\lab{lemma:20220928-zxbl1} Let $N=4$,
$b > \frac{1}{S^2}$ and $p \in (3,4)$. Then $c_1>0$.
\el
\bp
Let $A=\frac{a}{2}$, $B=\frac{1}{4}(b-\frac{1}{S^2})$, $\kappa_4=\frac{c^{4-p}}{2\|Q\|_{2}^{p-2}}$, $q=2(p-2)$ and $p_3\in (2,4)$. By \eqref{eq:20220926-e2}, we have that
\beq
\Omega_{A,B,0,\kappa_4}^{2,4,p_3,q}\rightarrow +\infty~\hbox{as}~c\rightarrow 0^+.
\eeq
Hence, by Lemma \ref{lemma:20220923-hhl1}, the continuity implies that there exists some $c_*>0$ such that $\displaystyle \Omega_{A,B,0,\kappa_4}^{2,4,p_3,q}> 1$ provided $c\leq c_*$. Then by the definition of $\displaystyle \Omega_{A,B,0,\kappa_4}^{2,4,p_3,q}$, the formula \eqref{eq:20220928-xe1} implies that $I(u)\geq 0$ for any $u\in S_c$ with $c\leq c_*$. And thus $I_c\geq 0, \forall c\in (0,c_*]$. Combining with Lemma \ref{lemma:20220928-zwl1}-(i), we obtain that $I_c=0, \forall c\in (0,c_*]$. Then by the definition of $c_1$, we have that $c_1\geq c_*>0$.
\ep

Similar to Theorem \ref{th:20220923-t1}, we can establish the following result for the normalized ground state solution.

\bt\lab{th:20220928-bt2} Let $N=4$,
$b > \frac{1}{S^2}$ and $p \in (3,4)$.
Then
\beq\lab{def:c1-another-2}
c_1=\sup\{c>0:I_c=0\}.
\eeq
And $I_c=0$ for $c \in (0,c_1]$ while $I_c<0$ for $c>c_1$.
In particular,
\begin{itemize}
\item[(i)]$I_c=0$ and it is not attained provided $0<c<c_1$.
\item[(ii)]$I_c<0$ and it is attained when $c> c_1$.
\item[(iii)] $I_{c_1}=0$ and it is attained.
\end{itemize}
\et

Similar to Theorem \ref{th:20220926-bt1} and Theorem \ref{th:20220926-bt2}, we get the following result for the normalized mountain pass type solution.

\bt\lab{th:20220928-bt1} Let $N=4$,
$b > \frac{1}{S^2}$ and $p \in (3,4)$. Then there exists $\bar{\eta}>0$ such that $I|_{S_{r,c}}$ has a mountain pass type critical point for $c\geq c_1-\bar{\eta}$.
\et

Similar to Theorem \ref{lemma:221225}, we get the following result about local constraint minimizer.

\bt\lab{th:1226-bt} Let $N=4$,
$b > \frac{1}{S^2}$ and $p \in (3,4)$. Then $I|_{S_{r,c}}$ has a local constraint minimizer with positive energy for $c \in  [c_1-\bar{\eta},c_1)$.
\et

We can also have the following non-existence result.

\bt\lab{th:20220928-bt3} Let $N=4$,
$b > \frac{1}{S^2}$ and $p \in (3,4)$. Then there exists $c_0>0$ such that problem \eqref{1.5}-\eqref{1.6} has no normalized solution provided $c< c_0$. In particular, one can check that
\begin{align}\label{eq:20220928-xbe1}
c_0\geq \frac{a \|Q\|_2^{\frac{p-2}{4-p}}}{(4-p)(p-2)^{\frac{1}{4-p}}}\left[\frac{1}{p-3}\left(b-\frac{1}{S^2}\right)\right]^{\frac{p-3}{4-p}}.
\end{align}
\et

\bp If $u\in H^1(\R^4)$ solves \eqref{1.5}-\eqref{1.6},
by (\ref{4.32}), we have that
\begin{align}\label{eq:20220928-xbe2}
a\|\nabla u\|_2^2+\left(b-\frac{1}{S^2}\right)\|\nabla u\|_2^4 \leq \frac{(p-2) c^{4-p}}{\|Q\|_2^{p-2}}\|\nabla u\|_2^{2(p-2)}.
\end{align}
By a similar argument in Lemma \ref{lemma:20220928-zxbl1}, there exists some $c_0>0$ such that
$$\Omega_{a, b-\frac{1}{S^2}, 0,\frac{(p-2) c^{4-p}}{\|Q\|_2^{p-2}}}^{2,4,p_3,2(p-2)}>1, \ \forall c< c_0.$$
Then by the definition of $\displaystyle \Omega_{a, b-\frac{1}{S^2}, 0,\frac{(p-2) c^{4-p}}{\|Q\|_2^{p-2}}}^{2,4,p_3,2(p-2)}$, the formula \eqref{eq:20220928-xbe2} implies that $\|\nabla u\|_2=0$ and thus $u\equiv 0$, a contradiction.
Furthermore, by \eqref{eq:20220926-e2}, a direct calculation shows that
\eqref{eq:20220928-xbe1} holds.
\ep

\noindent\textit{\bf Proof of  Theorem  \ref{th:main-zzzt2}-(i).} By Theorem \ref{th:20220928-bt2}, we get (i-1). By Theorem \ref{th:20220928-bt1} and Theorem \ref{th:1226-bt}, we get (i-2) and (i-3). By Theorem \ref{th:20220928-bt3}, we get (i-4).
\hfill $\square$

\subsection{The case of $p=3$}

\bl\lab{lemma:20220928-wl1}
Let $N=4, b>\frac{1}{S^2}$ and $p=3$. Then $c_1=a\|Q\|_2$.
\el
\bp
By \eqref{eq:20220928-xe1}, we have that
\beq\lab{eq:20220928-we1}
I(u)\geq \left(\frac{a}{2}-\frac{c}{2\|Q\|_2}\right)\|\nabla u\|_2^2 +\frac{1}{4}\left(b-\frac{1}{S^2}\right)\|\nabla u\|_2^4.
\eeq
Then for any $c<a\|Q\|_2$, the formula \eqref{eq:20220928-we1} implies that $I(u)>0$ for any $u\in S_c$ and thus $I_c\geq 0$. Combining with Lemma \ref{lemma:20220928-zwl1}-(i), we obtain that $I_c=0$. Hence, by the definition of $c_1$, we have that $c_1\geq a\|Q\|_2$.

On the other hand,  by \eqref{4.9}, we have that
\beq\lab{eq:20220928-we2}
I(Q_t)=\left(a-\frac{c}{\|Q\|_2}\right)\frac{1}{2}c^2t^2+\left(b-\frac{\|Q\|_4^4}{\|Q\|_2^4}\right)\frac{1}{4}c^4t^4.
\eeq
So for any $c>a\|Q||_2$, one can see that $I(Q_t)<0$ for $t>0$ small enough and thus $I_c\leq \min_{t>0}I(Q_t)<0$. So by the definition of $c_1$, we have that $c_1\leq a\|Q\|_2$.
\ep

\bt\lab{th:20220928-zhj-t1}
Let $N=4$, $b > \frac{1}{S^2}$ and $p=3$. Then $I_c<0$ is attained for $c>a\|Q\|_2$.
\et
\bp
By Lemma \ref{lemma:20220928-wl1}, we have that $I_c<0$ for $c>a\|Q\|_2$. Thanks to the compactness result in Lemma \ref{lemma:20220928-zwl3}, by a similar argument as Lemma \ref{lemma:20220923-xl1}, we can obtain the result.
\ep

Furthermore, we can prove the following non-existence result.
\bt\lab{th:20220928-zhj-t2} Let $N=4$,
$b > \frac{1}{S^2}$ and $p=3$. Then problem \eqref{1.5}-\eqref{1.6} has no normalized solution provided $c\leq  a\|Q\|_2$.
\et
\bp If $u\in H^1(\R^4)$ solves \eqref{1.5}-\eqref{1.6},
by \eqref{eq:20220928-xbe2}, we obtain that
\beq\lab{eq:20220928-we3}
\left(a-\frac{c}{\|Q\|_2}\right)\|\nabla u\|_2^2 +\left(b-\frac{1}{S^2}\right)\|\nabla u\|_2^4\leq 0.
\eeq
So for $c\leq a\|Q\|_2$, the formula \eqref{eq:20220928-we3} implies that $\|\nabla u\|_2=0$ and thus $u=0$, a contradiction.
\ep

\noindent\textit{\bf Proof of  Theorem  \ref{th:main-zzzt2}-(ii).} By Theorem \ref{th:20220928-zhj-t1} and Theorem \ref{th:20220928-zhj-t2}, we get the result.
\hfill $\square$

\subsection{The case of $2<p<3$}

\bl\lab{lemma:20220928-zbjs}
Let $N=4,b>\frac{1}{S^2}$ and $p \in (2,3)$. Then $I_c<0$ is attained for $c>0$.
\el
\bp
For any $c>0$, by \eqref{4.9}, we have that
\beq\lab{eq:20220928-bawe2}
I(Q_t)=\frac{a}{2} c^2 t^2 + \left(b-\frac{\|Q\|_4^4}{\|Q\|_2^4}\right)\frac{1}{4}c^4t^4-\frac{c^p t^{2p-4}}{2 \|Q\|_2^{p-2}}.
\eeq
Noting that $2p-4<2$, we have that $I(Q_t)<0$ for $t>0$ small enough and thus $I_c\leq \min_{t>0}I(Q_t)<0$. By a similar argument as Lemma \ref{lemma:20220923-xl1}, we can obtain that $I_c$ is attained by some $u\in S_{r,c}$.
\ep

\noindent\textit{\bf Proof of  Theorem  \ref{th:main-zzzt2}-(iii).} By Lemma \ref{lemma:20220928-zbjs}, we get the result.
\hfill $\square$

\end{document}